\documentclass[]{article}

\usepackage[margin=1in]{geometry}

\usepackage[usenames,dvipsnames]{xcolor}
\usepackage{wrapfig}
\usepackage{lettrine}
\usepackage{amssymb}
\usepackage{mdwlist}
\usepackage{color}
\usepackage{capt-of}
\usepackage{multirow}
\usepackage{graphicx}
\usepackage{epsfig}
\usepackage[euler]{textgreek}
\usepackage{epstopdf}
\usepackage{mathtools}
\usepackage{amsmath,amsthm,mathtools}
\usepackage{relsize}
\usepackage{breqn}
\usepackage{standalone}
\usepackage{tikz}
\usepackage{lineno}
\usepackage{textcomp}

\usepackage{caption}
\usepackage{subcaption}
\usepackage{soul}
\usepackage{cases}
\usepackage{afterpage}
\usepackage[section]{placeins}

\usepackage{hyperref}
\hypersetup{
    colorlinks,
    citecolor=blue,
    filecolor=black,
    linkcolor=blue,
    urlcolor=blue,
}
\usepackage{appendix}
\usepackage{color}

\usepackage{titling}
\usepackage[parfill]{parskip}
\usepackage{bm}

\numberwithin{equation}{section}

\newlength\tindent
\def\Snospace~{\S{}}

\allowdisplaybreaks[1]

\expandafter\def\expandafter\normalsize\expandafter{%
    \normalsize
    \setlength\abovedisplayskip{3pt}
    \setlength\belowdisplayskip{3pt}
    \setlength\abovedisplayshortskip{3pt}
    \setlength\belowdisplayshortskip{3pt}
}

\usepackage[T1]{fontenc}

\usepackage[breakable,theorems,skins]{tcolorbox}
\tcbset{enhanced}

\newcommand{\vect}[1]{\bm{{#1}}}
\newcommand{\tens}[1]{\underline{\underline{\bm{#1}}}}

\let\oldnabla\nabla
\renewcommand{\nabla}{{\vect{\oldnabla}}}

\usepackage[usenames,dvipsnames]{xcolor}
\newcommand{\makered}[1]{{\color{red}#1}}
\newcommand{\makegreen}[1]{{\color{ForestGreen}#1}}
\newcommand{\makeblue}[1]{{\color{blue}#1}}

 \title{On the Necessity of Superparametric Geometry Representation for Discontinuous Galerkin Methods on Domains with Curved Boundaries}

 \author{
  Philip Zwanenburg%
    \thanks{PhD Student. Email: philip.zwanenburg@mail.mcgill.ca}
   , and Siva Nadarajah%
    \thanks{Associate Professor, AIAA Associate Fellow. }
    \\
  {\normalsize\itshape Computational Aerodynamics Group}\\
  {\normalsize\itshape Department of Mechanical Engineering, McGill University}\\
  {\normalsize\itshape
  Montreal, QC H3A 2S6 Canada}
 }
 \date{May 3, 2017}

\begin{document}

\maketitle

\begin{abstract}
We provide numerical evidence demonstrating the necessity of employing a superparametric geometry representation in order to obtain optimal convergence orders on two-dimensional domains with curved boundaries when solving the Euler equations using Discontinuous Galerkin methods. However, concerning the obtention of optimal convergence orders for the Navier-Stokes equations, we demonstrate numerically that the use of isoparametric geometry representation is sufficient for the case considered here.

\end{abstract}

\section{Introduction}

\lettrine[nindent=0pt]{T}{he} proper treatment of complex geometry in high-order finite element methods has been shown to be crucial. That an isoparametric geometry representation is necessary for optimal error estimates when solving elliptic problems with Dirichlet or Neumann boundary conditions on curved domains using finite element methods has been rigorously established many years ago~\cite{lenoir1986,ciarlet1972}. Further, numerical results verifying optimal $O(h^{k+1})$ convergence of solution error in the $L^2$ norm demonstrate that this geometry representation is sufficient for a variety of equations and boundary conditions; the use of isoparametric finite elements in the context of Discontinuous Galerkin (DG) solvers has thus become a textbook approach~\cite{hesthaven2007}.

Despite this trend, it has been numerically observed that curved geometry must be represented using superparametric elements in the case of first order solution representation to maintain the optimal convergence properties of DG-type methods when solving the two-dimensional (2D) Euler equations~\cite{bassi1997, krivodonova2006}. As the superiority of high-order methods over industry-standard finite volume methods hinges on achieving high-order accurate convergence, several questions may be addressed following the preceding investigations. First, it is desired to determine whether the necessity of employing superparametric elements is restricted to the case of using a first order solution representation. It is also important to determine whether or not the same conclusions carry over when modelling the Euler equations in 3D or the Navier-Stokes equations.

Several of the questions above are addressed in this article. Notably, 2D numerical results presented below demonstrate that superparametric geometry representation is necessary for all solution orders when solving the Euler equations on domains with curved boundaries, notably when element aspect ratios deviate mildly from one. It is also demonstrated that isoparametric geometry representation is sufficient for the obtention of optimal convergence orders when solving the Navier-Stokes equations for the selected model problem. The major challenge which must be addressed if these results are to be extended to the 3D case would be to ensure sufficient mesh regularity on the refined mesh sequences. With this guarantee, it is expected that conclusions presented here would remain valid.

The paper is organized as follows. In~\autoref{sec:discretization}, we first discuss the DG discretization employed throughout this work. In~\autoref{sec:verification_Poisson}, verification of the implementation of the code when curved elements are considered is done using a model problem based on the Poisson equation. The first main contribution of the work is then presented in~\autoref{sec:verification_Euler} in which numerical results demonstrate that superparametric geometry representation is necessary for all solution orders when solving the Euler equations with curved slip wall boundary conditions. Finally, an analogous investigation is performed for the Navier-Stokes equations in~\autoref{sec:verification_NavierStokes} where it is demonstrated that isoparametric geometry representation is sufficient for the case considered here which employs adiabatic no-slip wall (Neumann) boundary conditions. Conclusions then follow in~\autoref{sec:conclusion}.

\section{Code \& Discretization Descriptions}
\label{sec:discretization}

All results presented in this work have been obtained with an in-house DG code freely available on~\href{https://github.com/PhilipZwanenburg/DPGSolver}{github}\footnote{https://github.com/PhilipZwanenburg/DPGSolver} and the corresponding author invites any questions related to their replication. The code allows for flexibility in the choice of element type (Tri, Quad, Tet, Hex, Wedge, Pyr), polynomial geometry representation order ($k_G$), as well as volume and face cubature order ($k_{vI}$ and $k_{fI}$). Unless otherwise indicated, the employed cubature rule was chosen such that polynomials of order $2k$ were integrated exactly. Following the standard DG methodology, the global solution is given by a composition of piecewise polynomial approximations where, on each element, the solution (conservative variables in the case of the Euler/Navier-Stokes equations) is represented by a combination of linearly independent modal or nodal polynomial basis functions of maximal order $k$; in the case of the Navier-Stokes equations we additionally employ the same order for the representation of the solution gradients. The geometry is represented in a globally $C^0$ polynomial space, whose order may be chosen independently of the solution order.

The implementation of the 2$^{nd}$ order equations (Poisson, Navier-Stokes) was done in a stabilized mixed formulation while the Euler equations were discretized in the standard weak form. Using the weak formulation precludes the need to represent the flux in a polynomial space when evaluating volume contributions to the residual; the flux is computed directly at the cubature nodes as required in the volume integral evaluation. This reduces aliasing errors as compared to discretizations based on the doubly integrated by parts weak form~\cite{hesthaven2007} in which aliasing errors are typically addressed through the use of the chain rule approach~\cite{wang2009}. 

As the test cases considered in this work all possess steady (time independent) solutions, an implicit solver with discretization based on the analytical linearization was used to converge the residuals to machine precision, using Petsc's~\cite{petsc2016} conjugate gradient (CG) and generalized minimal residual (GMRES) methods.

\section{Code Verification - Poisson Equation}
\label{sec:verification_Poisson}

\subsection{Discretization}

Consider the 3D Poisson equation
\begin{align} \label{eq:Poisson_std}
\nabla^2 u(\vect{x}) = s(\vect{x}),\ \vect{x} \coloneqq \begin{bmatrix} x & y & z \end{bmatrix} \in \vect{\Omega},
\end{align}

with source term, $s(\vect{x})$. Rewriting~\eqref{eq:Poisson_std} as a system of first order equations
\begin{align*}
\vect{q}(\vect{x}) = \nabla u(\vect{x}), \nabla \cdot \vect{q(\vect{x})} = s(\vect{x}),
\end{align*}

the two equations can then be integrated with respect to polynomial test functions, $\vect{\phi}_m(\vect{x})$, chosen as the basis functions, which are locally supported on each of the $m$ elements used to discretize the domain, to obtain the element-wise variational formulation of the problem,
\begin{align}
& \int_{\vect{\Omega}_v} \vect{\phi}_v(\vect{x})^T \vect{q}(\vect{x}) d \vect{\Omega}_v = \int_{\vect{\Omega}_v} \vect{\phi}_v(\vect{x})^T \nabla u(\vect{x}) d \vect{\Omega}_v \label{eq:discrete_Poisson_system1} \\
&  \int_{\vect{\Omega}_v} \vect{\phi}_v(\vect{x})^T  \nabla \cdot \vect{q(\vect{x})} d \vect{\Omega}_v =  \int_{\vect{\Omega}_v} \vect{\phi}_v(\vect{x})^T  s(\vect{x}) d \vect{\Omega}_v \label{eq:discrete_Poisson_system2}.
\end{align}

Performing integration by parts (IBP) twice in~\eqref{eq:discrete_Poisson_system1} and once
in~\eqref{eq:discrete_Poisson_system2}, transferring each element to the reference space and omitting the discrete solution, gradient and source dependence on $\vect{x}$ gives the mixed formulation
\begin{align}
& 
\int_{\vect{\Omega}_r} \vect{\phi}(\vect{r})^T J_v^{\Omega} \vect{q_h} d \vect{\Omega}_r = 
\int_{\vect{\Omega}_r} \vect{\phi}(\vect{r})^T J_v^{\Omega} \nabla u_h d \vect{\Omega}_r 
+ \int_{\vect{\Gamma}_r} \vect{\phi}(\vect{r})^T J_f^{\Gamma} \vect{\hat{n}}_f \left( u_h^{*} - u_h \right) d \vect{\Gamma}_r
 \label{eq:Poisson_strong_weak1} \\
&  
-\int_{\vect{\Omega}_r} J_v^{\Omega} \nabla \vect{\phi}(\vect{r})^T  \cdot \vect{q_h} d \vect{\Omega}_r 
+ \int_{\vect{\Gamma}_r} \vect{\phi}(\vect{r})^T J_f^{\Gamma} \vect{\hat{n}}_f \cdot \vect{q}_h^{*} d \vect{\Gamma}_r
=  \int_{\Omega_r} \vect{\phi}(\vect{r})^T  s_h d \vect{\Omega}_r
 \label{eq:Poisson_strong_weak2},
\end{align}

where $v$ and $f$ are the volume and face indices, $\vect{r} \coloneqq \begin{bmatrix} r & s & t \end{bmatrix}$ represents the coordinates in the reference space, star ($*$) superscripts are used to specify numerical traces and fluxes, and $ J_v^{\Omega}$ and $J_f^{\Gamma}$ denote the volume and area elements. The reader is invited to consult a previous publication~\cite{zwanenburg2016} for further details relating to the notation. The internal penalty method~\cite[Table \makeblue{1}]{brdar2012} was used for the computation of numerical traces and fluxes for the results presented in this section.

\subsection{Results}

The results presented here provide the verification of the code when treating curved elements. In all cases, meshes consisting entirely of triangular elements were used, with the solution, gradient and geometry represented in a nodal basis based on the alpha-optimized nodes~\cite{hesthaven2007} and the integration performed using the Gauss-Ledendre and Witherden-Vincent~\cite{witherden2015} nodes for face and volume integrals, respectively. Analogous results were obtained when using quadrilateral elements and are thus omitted for brevity.

The Poisson equation was solved on the quarter annulus with Dirichlet boundary conditions on the vertical and inner radial boundary and Neumann boundary conditions prescribed elsewhere with the source term chosen such that the exact solution is given by
\begin{align*}
u(x,y) = \sin(\pi x) \sin(\pi y).
\end{align*}

Results are presented for subparametric ($k_G = 1$), and isoparametric ($k_G = k$) geometry representation with the final geometry node placement defined according to an analytical mapping of a rectilinear domain to the quarter annulus; meshes generated in this way have the prefix ``ToBeCurved''. This can be contrasted with the more standard approach of only curving elements located adjacent to curved domain boundaries using a projection of face nodes to the boundary followed by a blending of internal nodes. The following conclusions can be drawn:

\begin{itemize}
\item Suboptimal convergence when using a subparametric geometry representation (\autoref{tab:Poisson_PG=1}).
\item Optimal convergence in the solution ($O(h^{k+1})$) and gradient ($O(h^{k})$) errors in $L^2$ when using an isoparametric geometry representation (\autoref{tab:Poisson_PG=P}).
\end{itemize}

\begin{table}[!h]
\begin{center}
\caption{ Errors and Convergence Orders - ToBeCurvedTRI meshes ($k_G = 1$)}
\label{tab:Poisson_PG=1}
\resizebox{\textwidth}{0.125\textheight}{
\begin{tabular}{| l | c | c c c | c c c | }
	\hline
	\multicolumn{2}{|c|}{} & \multicolumn{3}{c|}{$L^2$ Error} & \multicolumn{3}{c|}{Conv. Order} \\
	\hline
	Order ($k$) & Mesh Size ($h$) & $u$     & $q_1$   & $q_2$   & $u$     & $q_1$   & $q_2$   \\
	\hline
1	&  3.608e-02 &  3.844e-03 &  1.517e-01 &  1.862e-01 & - & - & - \\
	&  1.804e-02 &  9.800e-04 &  7.651e-02 &  9.379e-02 &  1.972e+00 &  9.876e-01 &  9.897e-01 \\
	&  9.021e-03 &  2.461e-04 &  3.838e-02 &  4.699e-02 &  1.994e+00 &  9.953e-01 &  9.970e-01 \\
	&  4.511e-03 &  6.158e-05 &  1.921e-02 &  2.351e-02 &  1.999e+00 &  9.983e-01 &  9.992e-01 \\
	\hline
2	&  2.552e-02 &  7.070e-04 &  6.611e-03 &  9.018e-03 & - & - & - \\
	&  1.276e-02 &  1.744e-04 &  1.808e-03 &  2.575e-03 &  2.020e+00 &  1.870e+00 &  1.808e+00 \\
	&  6.379e-03 &  4.336e-05 &  5.208e-04 &  7.704e-04 &  2.008e+00 &  1.796e+00 &  1.741e+00 \\
	&  3.189e-03 &  1.082e-05 &  1.590e-04 &  2.429e-04 &  2.003e+00 &  1.711e+00 &  1.665e+00 \\
	\hline
3	&  1.976e-02 &  7.006e-04 &  3.667e-03 &  3.769e-03 & - & - & - \\
	&  9.882e-03 &  1.736e-04 &  1.082e-03 &  1.216e-03 &  2.013e+00 &  1.761e+00 &  1.632e+00 \\
	&  4.941e-03 &  4.327e-05 &  3.396e-04 &  4.081e-04 &  2.004e+00 &  1.672e+00 &  1.574e+00 \\
	&  2.471e-03 &  1.080e-05 &  1.118e-04 &  1.403e-04 &  2.002e+00 &  1.603e+00 &  1.540e+00 \\
	\hline
\end{tabular}
}
\end{center}
\end{table}

\begin{table}[!h]
\begin{center}
\caption{ Errors and Convergence Orders - ToBeCurvedTRI meshes ($k_G = k$)}
\label{tab:Poisson_PG=P}
\resizebox{\textwidth}{0.125\textheight}{
\begin{tabular}{| l | c | c c c | c c c | }
	\hline
	\multicolumn{2}{|c|}{} & \multicolumn{3}{c|}{$L^2$ Error} & \multicolumn{3}{c|}{Conv. Order} \\
	\hline
	Order ($k$) & Mesh Size ($h$) & $u$     & $q_1$   & $q_2$   & $u$     & $q_1$   & $q_2$   \\
	\hline
1	&  3.608e-02 &  3.844e-03 &  1.517e-01 &  1.862e-01 & - & - & - \\
	&  1.804e-02 &  9.800e-04 &  7.651e-02 &  9.379e-02 &  1.972e+00 &  9.876e-01 &  9.897e-01 \\
	&  9.021e-03 &  2.461e-04 &  3.838e-02 &  4.699e-02 &  1.994e+00 &  9.953e-01 &  9.970e-01 \\
	&  4.511e-03 &  6.158e-05 &  1.921e-02 &  2.351e-02 &  1.999e+00 &  9.983e-01 &  9.992e-01 \\
	\hline
2	&  2.552e-02 &  8.871e-05 &  6.450e-03 &  8.334e-03 & - & - & - \\
	&  1.276e-02 &  1.132e-05 &  1.615e-03 &  2.134e-03 &  2.970e+00 &  1.998e+00 &  1.966e+00 \\
	&  6.379e-03 &  1.430e-06 &  4.039e-04 &  5.396e-04 &  2.985e+00 &  1.999e+00 &  1.983e+00 \\
	&  3.189e-03 &  1.796e-07 &  1.010e-04 &  1.356e-04 &  2.993e+00 &  2.000e+00 &  1.992e+00 \\
	\hline
3	&  1.976e-02 &  2.249e-06 &  2.074e-04 &  2.596e-04 & - & - & - \\
	&  9.882e-03 &  1.415e-07 &  2.605e-05 &  3.287e-05 &  3.991e+00 &  2.993e+00 &  2.981e+00 \\
	&  4.941e-03 &  8.853e-09 &  3.260e-06 &  4.135e-06 &  3.998e+00 &  2.999e+00 &  2.991e+00 \\
	&  2.471e-03 &  5.548e-10 &  4.075e-07 &  5.185e-07 &  3.996e+00 &  3.000e+00 &  2.996e+00 \\
	\hline
\end{tabular}
}
\end{center}
\end{table}

\section{On the Necessity of Superparametric Geometry Parametrization}

In this section we present results confirming the need for superparametric geometry representation when solving the Euler equations using slip-wall boundary conditions for a curved surface. We additionally present results demonstrating that an isoparametric geometry representation is sufficient when solving the Navier-Stokes equations using no-slip adiabatic wall boundary conditions for a curved surface for the specific case considered here.

\subsection{Discretization}

Following the notation of Toro~\cite[Chapter \makeblue{1}]{toro2009}, the steady 3D Navier-Stokes equations are given by
\begin{align} \label{eq:NavierStokes_std}
\nabla \cdot \left( \vect{F^i}(\vect{W}(\vect{x})) - \vect{F^v}(\vect{W}(\vect{x}),\vect{Q}(\vect{x})) \right) = \vect{0}^T, \vect{x} \coloneqq \begin{bmatrix} x & y & z \end{bmatrix} \in \vect{\Omega},
\end{align}

where the inviscid flux is defined as
\begin{align} \label{eq:F_inviscid}
\vect{F^i}(\vect{W}(\vect{x})) \coloneqq 
\begin{bmatrix} \rho\vect{v}^T & \rho \vect{v}^T \vect{v} + p \tens{I} & (E+p) \vect{v}^T \end{bmatrix},
\end{align}

and the viscous flux is defined as
\begin{align} \label{eq:F_viscous}
\vect{F^v}(\vect{W}(\vect{x}),\vect{Q}(\vect{x})) \coloneqq 
\begin{bmatrix} \vect{0}^T & \tens{\Pi} & \tens{\Pi} \vect{v}^T - \vect{q}^T \end{bmatrix}.
\end{align}

The various symbols represent the density, $\rho$, the velocity vector, $\vect{v}$, the total energy per unit volume, $E = \rho \left(e + \frac{1}{2} \vect{v} \vect{v}^T \right)$, the pressure, $p$, defined according to the equation of state for a calorically ideal gas,
\begin{align*}
p = (\gamma-1) \rho e,\ \gamma = \frac{C_p}{C_v}, 
\end{align*}

the stress tensor, $\tens{\Pi}$, given by
\begin{align*}
\tens{\Pi} = 2\mu \left( \tens{D} - \frac{1}{3} \nabla \cdot \vect{v} \tens{I} \right),\ \tens{D} \coloneqq \frac{1}{2} \left(\nabla^T \vect{v} + \left(\nabla^T \vect{v}\right)^T \right),
\end{align*}

where $\mu$ is the coefficient of shear viscosity and where the coefficient of bulk viscosity was assumed to be zero, and the energy flux vector, $\vect{q}$, defined by
\begin{align*}
\vect{q} = \kappa \nabla T,
\end{align*}

where $T$ represents the temperature and the coefficient of thermal conductivity is given by
\begin{align*}
\kappa = \frac{C_p \mu}{Pr},
\end{align*}

with $Pr$ representing the Prandtl number. In the case of the Euler equations, the contribution of the viscous flux is neglected. Defining a joint flux
\begin{align*}
\vect{F}(\vect{W}(\vect{x}),\vect{Q}(\vect{x})) \coloneqq 
\vect{F^i}(\vect{W}(\vect{x})) -
\vect{F^v}(\vect{W}(\vect{x}),\vect{Q}(\vect{x})),
\end{align*}

the discretization proceeds analogously to that for the case of the Poisson equation, performed in~\autoref{sec:verification_Poisson} such that one obtains,
\begin{align}
& 
\int_{\vect{\Omega}_r} \vect{\phi}(\vect{r})^T J_v^{\Omega} \vect{Q_h} d \vect{\Omega}_r = 
\int_{\vect{\Omega}_r} \vect{\phi}(\vect{r})^T J_v^{\Omega} \nabla \vect{W_h} d \vect{\Omega}_r 
+ \int_{\vect{\Gamma}_r} \vect{\phi}(\vect{r})^T J_f^{\Gamma} \vect{\hat{n}}_f \left( \vect{W_h}^{*} - \vect{W_h} \right) d \vect{\Gamma}_r
 \label{eq:NS_MF1} \\
&  
-\int_{\vect{\Omega}_r} J_v^{\Omega} \nabla \vect{\phi}(\vect{r})^T  \cdot \vect{F_h} d \vect{\Omega}_r 
+ \int_{\vect{\Gamma}_r} \vect{\phi}(\vect{r})^T J_f^{\Gamma} \vect{\hat{n}}_f \cdot \vect{F}_h^{*} d \vect{\Gamma}_r
=  \vect{0}^T
 \label{eq:NS_MF2}.
\end{align}

In the case of the Euler equations,~\eqref{eq:NS_MF1} is neglected, and we set $\vect{F} = \vect{F^i}$. For the results presented below, the inviscid numerical fluxes were computed using the Roe-Pike method~\cite[Chapter \makeblue{11.3}]{toro2009} and the viscous numerical traces and fluxes were computed using the BR2 method~\cite[eq. (\makeblue{10})]{bassi2010},~\cite[eq. (\makeblue{4.3})]{brdar2012} with the penalty scaling parameter set equal to the maximum number of element faces, such that the coercivity estimate holds~\cite{brezzi2000}.

\subsection{Results - Euler}
\label{sec:verification_Euler}

For the verification of the main result, a 2D inviscid isentropic supersonic flow in a quarter annulus with inner radius, $r_{i} = 1$, and outer radius, $r_{o} = 1.384$, in the first quadrant was considered~\cite{krivodonova2006}. The density and Mach number on the inner surface are set to $\rho_{i} = 1$ and $M_{i} = 2.25$, with the exact density, pressure and velocity magnitude distributions given by
\begin{align*}
\rho = \rho_{i} \left( 1 + \frac{\gamma-1}{2} M_{i}^2 \left( 1 - \left( \frac{r_{i}}{r} \right)^2 \right) \right)^{\frac{1}{\gamma-1}},\
p = \frac{\rho^{\gamma}}{\gamma},\
|| \vect{v} || = \frac{M_{i}}{r}.
\end{align*}

Riemann invariant boundary conditions were applied at the inlet and outlet (equivalent to supersonic inflow/outflow boundary conditions in this case) and slip wall boundary conditions were applied along the curved walls. As the analytical solution is available, the integrated $L^2$ errors of all primitive variables as well as of the normalized entropy were measured, with the entropy error given by
\begin{align*}
|| \Delta s ||_{L^2(\vect{\Omega})} = \left( \frac{1}{| \vect{\Omega} |} \int_{\vect{\Omega}} \left(\frac{\frac{p}{\rho^{\gamma}}-\frac{p_0}{{\rho_0}^{\gamma}}}{\frac{p_0}{{\rho_0}^{\gamma}}} \right)^2 d \vect{\Omega} \right)^{\frac{1}{2}}.
\end{align*}

\begin{figure}[!h]
    \centering
    \begin{subfigure}{0.485\textwidth}
        \includegraphics[width=\textwidth]{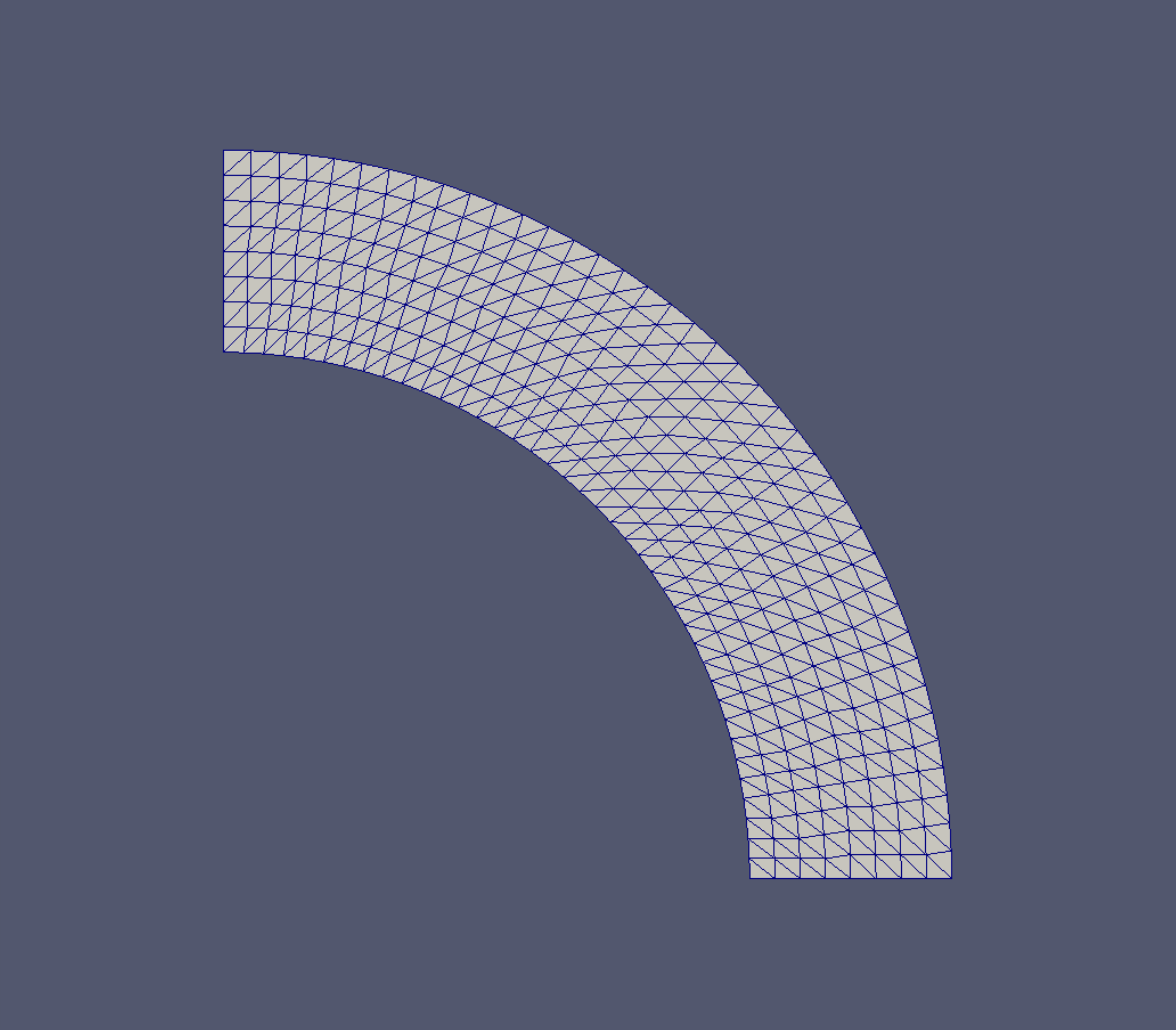}
        \caption{AR $ \approx 1.0$}
    \end{subfigure}
\hfill
    \begin{subfigure}{0.485\textwidth}
        \includegraphics[width=\textwidth]{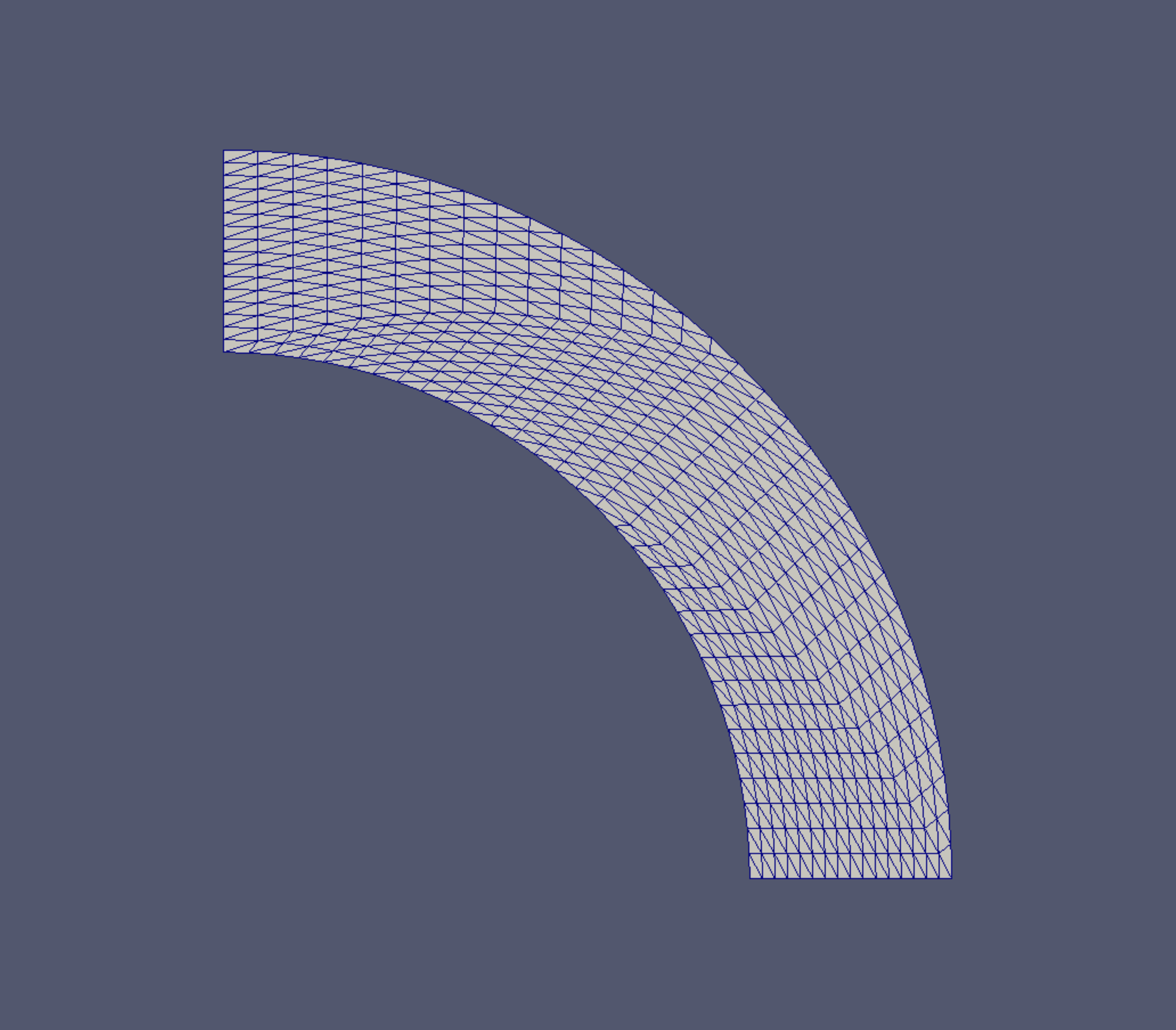}
        \caption{AR $ \approx 2.5$}
    \end{subfigure}
\vskip\baselineskip
    \begin{subfigure}{0.485\textwidth}
        \includegraphics[width=\textwidth]{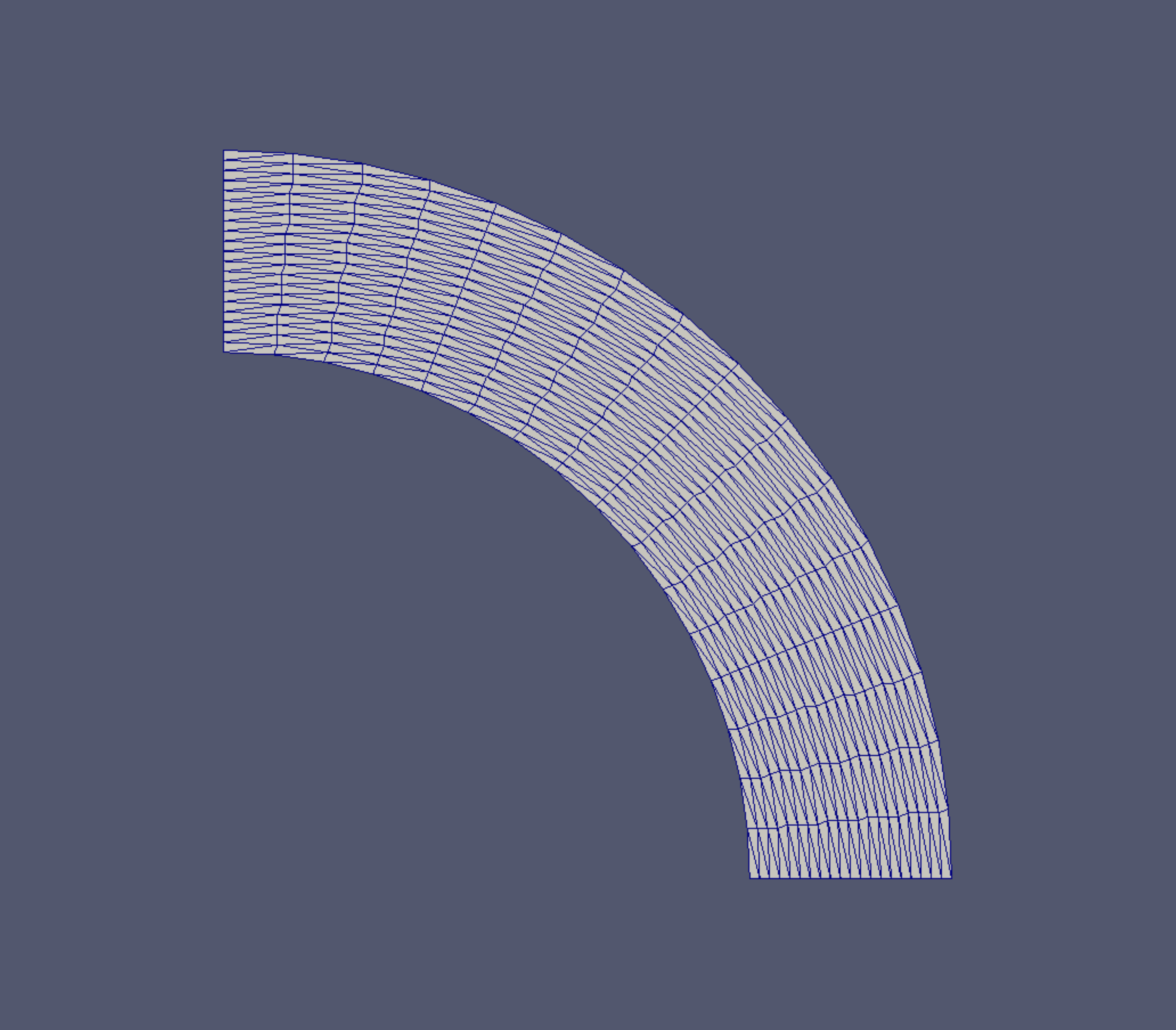}
        \caption{AR $ \approx 5.0$}
    \end{subfigure}
\hfill
    \begin{subfigure}{0.485\textwidth}
        \includegraphics[width=\textwidth]{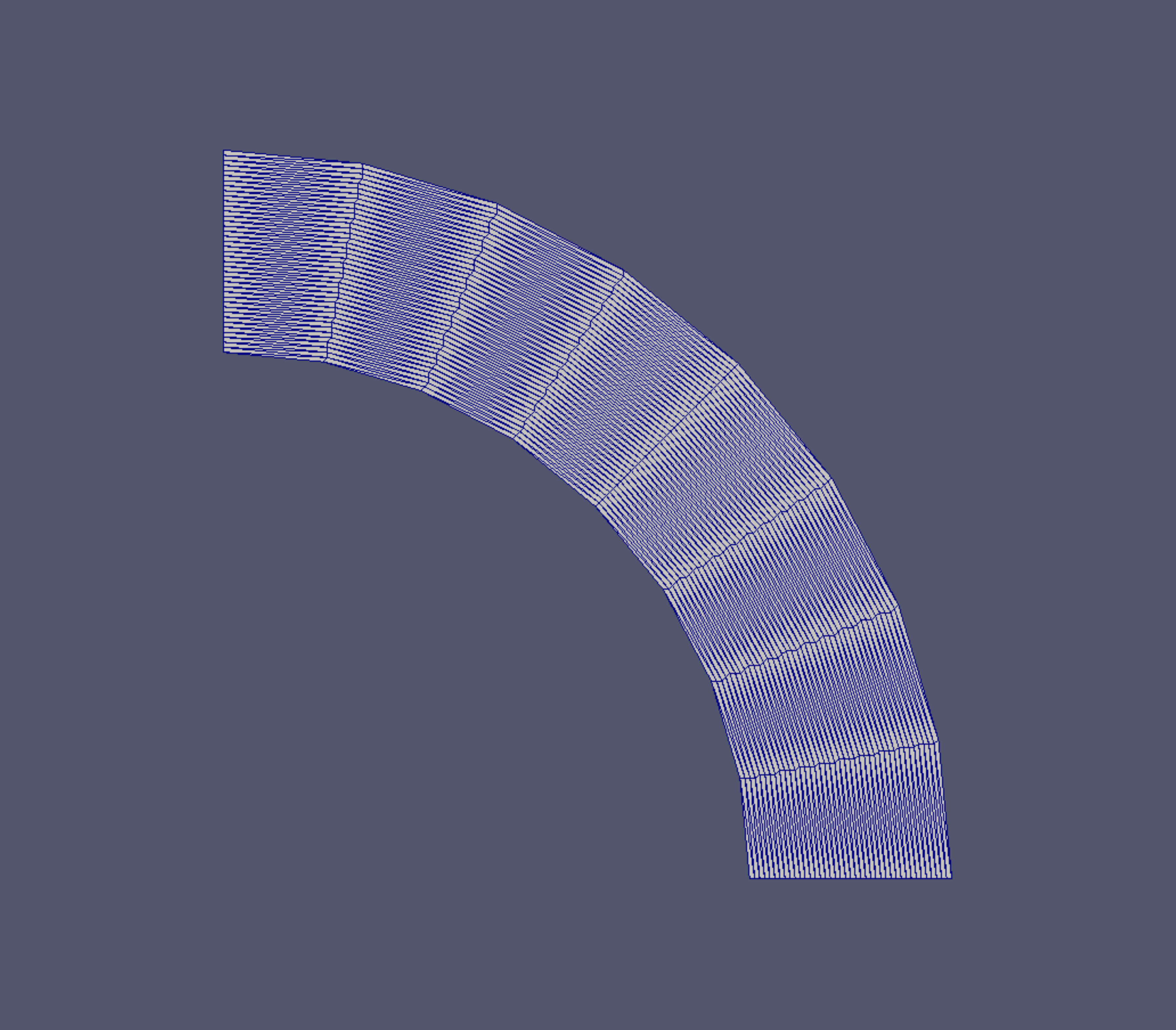}
        \caption{AR $ \approx 20.0$}
    \end{subfigure}
    \caption{Sample Meshes used in the Euler Convergence Studies}
    \label{fig:SupersonicVortex_Meshes}
\end{figure}

We present results on meshes with element aspect ratios (AR) ranging approximately from one to twenty to establish our conclusions; sample meshes can be seen in~\autoref{fig:SupersonicVortex_Meshes}.

The following conclusions can be drawn:

\begin{itemize}
\item The previous conclusion~\cite{bassi1997} claiming that superparametric geometry representation is only necessary for $k1$ solution representation is valid if only entropy convergence is taken into account and elements with aspect ratio close to one are used (\autoref{tab:Euler_PG=P_AR=1.0} and~\autoref{tab:Euler_PG=P_AR=2.5}). 
\item The sensitivity of the entropy error convergence to the geometry representation is less than that of the primitive variables (compare~\autoref{tab:Euler_PG=P_AR=2.5} with~\autoref{tab:Euler_PG=P_AR=5.0} for $k > 1$).
\item That optimal convergence orders are still obtained when extremely large aspect ratio elements are present when the superparametric geometry is used indicates that $k_G = k+1$ is sufficient (\autoref{tab:Euler_PG=P+1_AR=20.0}).
\end{itemize}

Taken together, the conclusions above suggest that extremely high quality meshes should not be used for convergence order studies, and furthermore, that entropy is not as good of an indicator of correct curved element treatment as the primitive variables, but may still be used  if meshes formed from elements with high enough aspect ratios (AR greater than 2.5 in this case) are used.

\begin{table}[!h]
\begin{center}
\caption{ Errors and Convergence Orders - ToBeCurvedTRI meshes ($k_G = k$, AR $\approx$ 1.0)}
\label{tab:Euler_PG=P_AR=1.0}
\resizebox{\textwidth}{!}{
\begin{tabular}{| l | c | c c c c c | c c c c c | }
	\hline
\multicolumn{2}{|c|}{} & \multicolumn{5}{c|}{$L^2$ Error} & \multicolumn{5}{c|}{Conv. Order} \\
	\hline
	Order ($k$) & Mesh Size ($h$) & $\rho$ & $u$     & $v$     & $p$     & $s$     & $\rho$ & $u$     & $v$     & $p$     & $s$     \\
	\hline
1	&  4.564e-02 &  1.516e-02 &  1.143e-02 &  1.338e-02 &  2.281e-02 &  5.630e-03 & - & - & - & - & - \\
	&  2.282e-02 &  4.506e-03 &  3.674e-03 &  4.101e-03 &  6.544e-03 &  1.965e-03 &  1.750e+00 &  1.637e+00 &  1.705e+00 &  1.801e+00 &  1.519e+00 \\
	&  1.141e-02 &  1.449e-03 &  1.234e-03 &  1.335e-03 &  2.028e-03 &  6.988e-04 &  1.637e+00 &  1.574e+00 &  1.620e+00 &  1.690e+00 &  1.491e+00 \\
	&  5.705e-03 &  4.874e-04 &  4.235e-04 &  4.484e-04 &  6.616e-04 &  2.494e-04 &  \makered{1.572e+00} &  \makered{1.543e+00} &  \makered{1.574e+00} &  \makered{1.616e+00} &  \makered{1.486e+00} \\
	\hline
2	&  3.227e-02 &  1.918e-04 &  6.501e-05 &  1.681e-04 &  1.092e-04 &  1.951e-04 & - & - & - & - & - \\
	&  1.614e-02 &  2.565e-05 &  9.413e-06 &  2.236e-05 &  1.460e-05 &  2.639e-05 &  2.902e+00 &  2.788e+00 &  2.910e+00 &  2.903e+00 &  2.886e+00 \\
	&  8.069e-03 &  3.394e-06 &  1.439e-06 &  2.894e-06 &  2.164e-06 &  3.401e-06 &  2.918e+00 &  2.709e+00 &  2.950e+00 &  2.754e+00 &  2.956e+00 \\
	&  4.034e-03 &  4.548e-07 &  2.319e-07 &  3.737e-07 &  3.435e-07 &  4.302e-07 &  \makegreen{2.900e+00} &  \makered{2.634e+00} &  \makegreen{2.953e+00} &  \makered{2.655e+00} &  \makegreen{2.983e+00} \\
	\hline
3	&  2.500e-02 &  4.057e-06 &  1.491e-06 &  3.779e-06 &  2.164e-06 &  4.254e-06 & - & - & - & - & - \\
	&  1.250e-02 &  2.622e-07 &  9.803e-08 &  2.426e-07 &  1.288e-07 &  2.789e-07 &  3.952e+00 &  3.927e+00 &  3.962e+00 &  4.071e+00 &  3.931e+00 \\
	&  6.250e-03 &  1.687e-08 &  6.360e-09 &  1.550e-08 &  8.176e-09 &  1.795e-08 &  3.959e+00 &  3.946e+00 &  3.968e+00 &  3.977e+00 &  3.958e+00 \\
	&  3.125e-03 &  1.074e-09 &  4.107e-10 &  9.832e-10 &  5.270e-10 &  1.140e-09 &  \makegreen{3.973e+00} &  \makegreen{3.953e+00} &  \makegreen{3.979e+00} &  \makegreen{3.956e+00} &  \makegreen{3.977e+00} \\
	\hline
\end{tabular}
}
\end{center}
\end{table}

\begin{table}[!h]
\begin{center}
\caption{ Errors and Convergence Orders - ToBeCurvedTRI meshes ($k_G = k+1$, AR $\approx$ 1.0)}
\label{tab:Euler_PG=P+1_AR=1.0}
\resizebox{\textwidth}{!}{
\begin{tabular}{| l | c | c c c c c | c c c c c | }
	\hline
\multicolumn{2}{|c|}{} & \multicolumn{5}{c|}{$L^2$ Error} & \multicolumn{5}{c|}{Conv. Order} \\
	\hline
	Order ($k$) & Mesh Size ($h$) & $\rho$ & $u$     & $v$     & $p$     & $s$     & $\rho$ & $u$     & $v$     & $p$     & $s$     \\
	\hline
1	&  4.564e-02 &  2.869e-03 &  1.090e-03 &  2.510e-03 &  1.484e-03 &  3.064e-03 & - & - & - & - & - \\
	&  2.282e-02 &  7.932e-04 &  2.898e-04 &  7.070e-04 &  3.658e-04 &  8.699e-04 &  1.855e+00 &  1.911e+00 &  1.828e+00 &  2.021e+00 &  1.817e+00 \\
	&  1.141e-02 &  2.076e-04 &  7.446e-05 &  1.870e-04 &  8.898e-05 &  2.304e-04 &  1.934e+00 &  1.960e+00 &  1.919e+00 &  2.040e+00 &  1.917e+00 \\
	&  5.705e-03 &  5.317e-05 &  1.886e-05 &  4.814e-05 &  2.187e-05 &  5.931e-05 &  1.965e+00 &  1.981e+00 &  1.958e+00 &  2.025e+00 &  1.958e+00 \\
	\hline
2	&  3.227e-02 &  1.811e-04 &  4.811e-05 &  1.564e-04 &  9.082e-05 &  1.852e-04 & - & - & - & - & - \\
	&  1.614e-02 &  2.385e-05 &  6.092e-06 &  2.078e-05 &  1.063e-05 &  2.510e-05 &  2.924e+00 &  2.981e+00 &  2.912e+00 &  3.095e+00 &  2.883e+00 \\
	&  8.069e-03 &  3.060e-06 &  7.669e-07 &  2.663e-06 &  1.308e-06 &  3.235e-06 &  2.963e+00 &  2.990e+00 &  2.964e+00 &  3.023e+00 &  2.956e+00 \\
	&  4.034e-03 &  3.868e-07 &  9.591e-08 &  3.361e-07 &  1.625e-07 &  4.090e-07 &  2.984e+00 &  2.999e+00 &  2.986e+00 &  3.009e+00 &  2.984e+00 \\
	\hline
3	&  2.500e-02 &  3.919e-06 &  1.429e-06 &  3.554e-06 &  2.139e-06 &  4.026e-06 & - & - & - & - & - \\
	&  1.250e-02 &  2.525e-07 &  9.343e-08 &  2.279e-07 &  1.265e-07 &  2.639e-07 &  3.956e+00 &  3.935e+00 &  3.963e+00 &  4.080e+00 &  3.931e+00 \\
	&  6.250e-03 &  1.620e-08 &  5.997e-09 &  1.457e-08 &  7.923e-09 &  1.699e-08 &  3.962e+00 &  3.962e+00 &  3.968e+00 &  3.997e+00 &  3.957e+00 \\
	&  3.125e-03 &  1.027e-09 &  3.787e-10 &  9.216e-10 &  4.973e-10 &  1.078e-09 &  3.979e+00 &  3.985e+00 &  3.982e+00 &  3.994e+00 &  3.978e+00 \\
	\hline
\end{tabular}
}
\end{center}
\end{table}

\begin{table}[!h]
\begin{center}
\caption{ Errors and Convergence Orders - ToBeCurvedTRI meshes ($k_G = k$, AR $\approx$ 2.5)}
\label{tab:Euler_PG=P_AR=2.5}
\resizebox{\textwidth}{!}{
\begin{tabular}{| l | c | c c c c c | c c c c c | }
	\hline
\multicolumn{2}{|c|}{} & \multicolumn{5}{c|}{$L^2$ Error} & \multicolumn{5}{c|}{Conv. Order} \\
	\hline
	Order ($k$) & Mesh Size ($h$) & $\rho$ & $u$     & $v$     & $p$     & $s$     & $\rho$ & $u$     & $v$     & $p$     & $s$     \\
	\hline
1	&  3.608e-02 &  2.829e-02 &  2.232e-02 &  2.619e-02 &  4.254e-02 &  1.245e-02 & - & - & - & - & - \\
	&  1.804e-02 &  9.366e-03 &  7.693e-03 &  8.306e-03 &  1.328e-02 &  4.217e-03 &  1.595e+00 &  1.537e+00 &  1.657e+00 &  1.680e+00 &  1.562e+00 \\
	&  9.021e-03 &  3.291e-03 &  2.754e-03 &  2.817e-03 &  4.526e-03 &  1.452e-03 &  1.509e+00 &  1.482e+00 &  1.560e+00 &  1.552e+00 &  1.538e+00 \\
	&  4.511e-03 &  1.179e-03 &  9.964e-04 &  9.814e-04 &  1.599e-03 &  5.056e-04 &  1.481e+00 &  1.467e+00 &  1.521e+00 &  1.501e+00 &  1.522e+00 \\
	\hline
2	&  2.552e-02 &  1.622e-04 &  1.310e-04 &  1.206e-04 &  2.193e-04 &  4.701e-05 & - & - & - & - & - \\
	&  1.276e-02 &  3.005e-05 &  2.405e-05 &  2.241e-05 &  4.133e-05 &  6.266e-06 &  2.432e+00 &  2.445e+00 &  2.428e+00 &  2.408e+00 &  2.907e+00 \\
	&  6.379e-03 &  5.607e-06 &  4.477e-06 &  4.175e-06 &  7.783e-06 &  8.190e-07 &  2.422e+00 &  2.426e+00 &  2.424e+00 &  2.409e+00 &  2.936e+00 \\
	&  3.189e-03 &  1.041e-06 &  8.336e-07 &  7.711e-07 &  1.451e-06 &  1.066e-07 &  \makered{2.429e+00} &  \makered{2.425e+00} &  \makered{2.437e+00} &  \makered{2.423e+00} &  \makegreen{2.941e+00} \\
	\hline
3	&  1.976e-02 &  1.052e-06 &  8.885e-07 &  7.834e-07 &  1.383e-06 &  4.381e-07 & - & - & - & - & - \\
	&  9.882e-03 &  9.655e-08 &  7.681e-08 &  7.440e-08 &  1.308e-07 &  2.902e-08 &  3.445e+00 &  3.532e+00 &  3.396e+00 &  3.402e+00 &  3.916e+00 \\
	&  4.941e-03 &  9.732e-09 &  7.669e-09 &  7.498e-09 &  1.344e-08 &  1.889e-09 &  3.310e+00 &  3.324e+00 &  3.311e+00 &  3.283e+00 &  3.941e+00 \\
	&  2.471e-03 &  9.769e-10 &  7.525e-10 &  7.586e-10 &  1.360e-09 &  1.227e-10 &  \makered{3.316e+00} &  \makered{3.349e+00} &  \makered{3.305e+00} &  \makered{3.305e+00} &  \makegreen{3.945e+00} \\
	\hline
\end{tabular}
}
\end{center}
\end{table}

\begin{table}[!h]
\begin{center}
\caption{ Errors and Convergence Orders - ToBeCurvedTRI meshes ($k_G = k+1$, AR $\approx$ 2.5)}
\label{tab:Euler_PG=P+1_AR=2.5}
\resizebox{\textwidth}{!}{
\begin{tabular}{| l | c | c c c c c | c c c c c | }
	\hline
\multicolumn{2}{|c|}{} & \multicolumn{5}{c|}{$L^2$ Error} & \multicolumn{5}{c|}{Conv. Order} \\
	\hline
	Order ($k$) & Mesh Size ($h$) & $\rho$ & $u$     & $v$     & $p$     & $s$     & $\rho$ & $u$     & $v$     & $p$     & $s$     \\
	\hline
1	&  3.608e-02 &  1.746e-03 &  8.055e-04 &  1.394e-03 &  1.185e-03 &  1.805e-03 & - & - & - & - & - \\
	&  1.804e-02 &  4.688e-04 &  2.095e-04 &  3.789e-04 &  2.815e-04 &  4.997e-04 &  1.897e+00 &  1.943e+00 &  1.880e+00 &  2.074e+00 &  1.853e+00 \\
	&  9.021e-03 &  1.226e-04 &  5.349e-05 &  9.958e-05 &  6.825e-05 &  1.325e-04 &  1.935e+00 &  1.970e+00 &  1.928e+00 &  2.044e+00 &  1.915e+00 \\
	&  4.511e-03 &  3.149e-05 &  1.353e-05 &  2.562e-05 &  1.682e-05 &  3.426e-05 &  1.962e+00 &  1.983e+00 &  1.959e+00 &  2.021e+00 &  1.952e+00 \\
	\hline
2	&  2.552e-02 &  4.723e-05 &  1.767e-05 &  3.924e-05 &  2.896e-05 &  4.636e-05 & - & - & - & - & - \\
	&  1.276e-02 &  6.104e-06 &  2.243e-06 &  5.088e-06 &  3.496e-06 &  6.116e-06 &  2.952e+00 &  2.977e+00 &  2.947e+00 &  3.050e+00 &  2.922e+00 \\
	&  6.379e-03 &  7.756e-07 &  2.801e-07 &  6.472e-07 &  4.209e-07 &  7.869e-07 &  2.976e+00 &  3.002e+00 &  2.975e+00 &  3.054e+00 &  2.958e+00 \\
	&  3.189e-03 &  9.777e-08 &  3.496e-08 &  8.161e-08 &  5.138e-08 &  9.982e-08 &  2.988e+00 &  3.002e+00 &  2.987e+00 &  3.034e+00 &  2.979e+00 \\
	\hline
3	&  1.976e-02 &  4.794e-07 &  2.445e-07 &  4.494e-07 &  4.268e-07 &  4.433e-07 & - & - & - & - & - \\
	&  9.882e-03 &  3.223e-08 &  1.619e-08 &  3.017e-08 &  2.944e-08 &  2.926e-08 &  3.895e+00 &  3.916e+00 &  3.897e+00 &  3.858e+00 &  3.921e+00 \\
	&  4.941e-03 &  2.061e-09 &  1.035e-09 &  1.916e-09 &  1.869e-09 &  1.875e-09 &  3.967e+00 &  3.968e+00 &  3.977e+00 &  3.977e+00 &  3.964e+00 \\
	&  2.471e-03 &  1.267e-10 &  6.407e-11 &  1.173e-10 &  1.114e-10 &  1.178e-10 &  \makegreen{4.024e+00} &  \makegreen{4.014e+00} &  \makegreen{4.030e+00} &  \makegreen{4.068e+00} &  \makegreen{3.992e+00} \\
	\hline
\end{tabular}
}
\end{center}
\end{table}

\begin{table}[!h]
\begin{center}
\caption{ Errors and Convergence Orders - ToBeCurvedTRI meshes ($k_G = k$, AR $\approx$ 5.0)}
\label{tab:Euler_PG=P_AR=5.0}
\resizebox{\textwidth}{!}{
\begin{tabular}{| l | c | c c c c c | c c c c c | }
	\hline
\multicolumn{2}{|c|}{} & \multicolumn{5}{c|}{$L^2$ Error} & \multicolumn{5}{c|}{Conv. Order} \\
	\hline
	Order ($k$) & Mesh Size ($h$) & $\rho$ & $u$     & $v$     & $p$     & $s$     & $\rho$ & $u$     & $v$     & $p$     & $s$     \\
	\hline
1	&  4.564e-02 &  9.000e-02 &  6.665e-02 &  1.011e-01 &  1.484e-01 &  5.212e-02 & - & - & - & - & - \\
	&  2.282e-02 &  2.542e-02 &  2.053e-02 &  2.826e-02 &  3.646e-02 &  1.874e-02 &  1.824e+00 &  1.699e+00 &  1.839e+00 &  2.025e+00 &  1.476e+00 \\
	&  1.141e-02 &  8.253e-03 &  6.823e-03 &  9.057e-03 &  1.065e-02 &  6.683e-03 &  1.623e+00 &  1.589e+00 &  1.642e+00 &  1.776e+00 &  1.488e+00 \\
	&  5.705e-03 &  2.808e-03 &  2.344e-03 &  3.046e-03 &  3.368e-03 &  2.371e-03 &  1.556e+00 &  1.541e+00 &  1.572e+00 &  1.660e+00 &  1.495e+00 \\
	\hline
2	&  3.227e-02 &  6.941e-04 &  6.442e-04 &  4.495e-04 &  9.731e-04 &  8.156e-05 & - & - & - & - & - \\
	&  1.614e-02 &  1.249e-04 &  1.135e-04 &  8.300e-05 &  1.748e-04 &  1.128e-05 &  2.474e+00 &  2.505e+00 &  2.437e+00 &  2.477e+00 &  2.854e+00 \\
	&  8.069e-03 &  2.205e-05 &  1.984e-05 &  1.474e-05 &  3.085e-05 &  1.634e-06 &  2.502e+00 &  2.516e+00 &  2.493e+00 &  2.503e+00 &  2.787e+00 \\
	&  4.034e-03 &  3.893e-06 &  3.497e-06 &  2.599e-06 &  5.446e-06 &  2.481e-07 &  2.502e+00 &  2.504e+00 &  2.504e+00 &  2.502e+00 &  \makered{2.719e+00} \\
	\hline
3	&  2.500e-02 &  9.011e-06 &  8.724e-06 &  5.790e-06 &  1.263e-05 &  1.295e-06 & - & - & - & - & - \\
	&  1.250e-02 &  9.565e-07 &  8.391e-07 &  6.944e-07 &  1.336e-06 &  1.134e-07 &  3.236e+00 &  3.378e+00 &  3.060e+00 &  3.241e+00 &  3.514e+00 \\
	&  6.250e-03 &  8.914e-08 &  7.223e-08 &  6.850e-08 &  1.245e-07 &  9.393e-09 &  3.424e+00 &  3.538e+00 &  3.342e+00 &  3.424e+00 &  3.594e+00 \\
	&  3.125e-03 &  8.437e-09 &  6.640e-09 &  6.523e-09 &  1.178e-08 &  8.085e-10 &  3.401e+00 &  3.443e+00 &  3.392e+00 &  3.401e+00 &  \makered{3.538e+00} \\
	\hline
\end{tabular}
}
\end{center}
\end{table}

\begin{table}[!h]
\begin{center}
\caption{ Errors and Convergence Orders - ToBeCurvedTRI meshes ($k_G = k+1$, AR $\approx$ 5.0)}
\label{tab:Euler_PG=P+1_AR=5.0}
\resizebox{\textwidth}{!}{
\begin{tabular}{| l | c | c c c c c | c c c c c | }
	\hline
\multicolumn{2}{|c|}{} & \multicolumn{5}{c|}{$L^2$ Error} & \multicolumn{5}{c|}{Conv. Order} \\
	\hline
	Order ($k$) & Mesh Size ($h$) & $\rho$ & $u$     & $v$     & $p$     & $s$     & $\rho$ & $u$     & $v$     & $p$     & $s$     \\
	\hline
1	&  4.564e-02 &  2.838e-03 &  2.225e-03 &  2.479e-03 &  3.167e-03 &  3.199e-03 & - & - & - & - & - \\
	&  2.282e-02 &  7.800e-04 &  5.847e-04 &  7.034e-04 &  7.936e-04 &  9.069e-04 &  1.863e+00 &  1.928e+00 &  1.817e+00 &  1.997e+00 &  1.818e+00 \\
	&  1.141e-02 &  2.069e-04 &  1.505e-04 &  1.892e-04 &  1.982e-04 &  2.437e-04 &  1.915e+00 &  1.958e+00 &  1.895e+00 &  2.002e+00 &  1.896e+00 \\
	&  5.705e-03 &  5.347e-05 &  3.823e-05 &  4.919e-05 &  4.943e-05 &  6.340e-05 &  1.952e+00 &  1.977e+00 &  1.943e+00 &  2.003e+00 &  1.943e+00 \\
	\hline
2	&  3.227e-02 &  6.915e-05 &  4.802e-05 &  7.681e-05 &  4.366e-05 &  7.775e-05 & - & - & - & - & - \\
	&  1.614e-02 &  8.960e-06 &  6.112e-06 &  1.018e-05 &  5.048e-06 &  1.032e-05 &  2.948e+00 &  2.974e+00 &  2.915e+00 &  3.113e+00 &  2.913e+00 \\
	&  8.069e-03 &  1.141e-06 &  7.710e-07 &  1.312e-06 &  5.960e-07 &  1.331e-06 &  2.973e+00 &  2.987e+00 &  2.956e+00 &  3.082e+00 &  2.955e+00 \\
	&  4.034e-03 &  1.440e-07 &  9.684e-08 &  1.665e-07 &  7.193e-08 &  1.691e-07 &  2.987e+00 &  2.993e+00 &  2.978e+00 &  3.051e+00 &  2.977e+00 \\
	\hline
3	&  2.500e-02 &  7.920e-07 &  5.872e-07 &  9.246e-07 &  6.707e-07 &  8.109e-07 & - & - & - & - & - \\
	&  1.250e-02 &  5.204e-08 &  3.884e-08 &  6.086e-08 &  4.267e-08 &  5.362e-08 &  3.928e+00 &  3.919e+00 &  3.925e+00 &  3.974e+00 &  3.919e+00 \\
	&  6.250e-03 &  3.284e-09 &  2.448e-09 &  3.896e-09 &  2.561e-09 &  3.450e-09 &  3.986e+00 &  3.987e+00 &  3.965e+00 &  4.059e+00 &  3.958e+00 \\
	&  3.125e-03 &  2.033e-10 &  1.537e-10 &  2.438e-10 &  1.491e-10 &  2.186e-10 &  4.014e+00 &  3.994e+00 &  3.999e+00 &  4.102e+00 &  3.980e+00 \\
	\hline
\end{tabular}
}
\end{center}
\end{table}

\begin{table}[!h]
\begin{center}
\caption{ Errors and Convergence Orders - ToBeCurvedTRI meshes ($k_G = k+1$, AR $\approx$ 20.0)}
\label{tab:Euler_PG=P+1_AR=20.0}
\resizebox{\textwidth}{!}{
\begin{tabular}{| l | c | c c c c c | c c c c c | }
	\hline
\multicolumn{2}{|c|}{} & \multicolumn{5}{c|}{$L^2$ Error} & \multicolumn{5}{c|}{Conv. Order} \\
	\hline
	Order ($k$) & Mesh Size ($h$) & $\rho$ & $u$     & $v$     & $p$     & $s$     & $\rho$ & $u$     & $v$     & $p$     & $s$     \\
	\hline
1	&  4.564e-02 &  7.462e-03 &  8.536e-03 &  7.392e-03 &  1.335e-02 &  7.989e-03 & - & - & - & - & - \\
	&  2.282e-02 &  2.061e-03 &  2.002e-03 &  1.933e-03 &  3.315e-03 &  2.257e-03 &  1.856e+00 &  2.092e+00 &  1.935e+00 &  2.010e+00 &  1.824e+00 \\
	&  1.141e-02 &  5.652e-04 &  5.083e-04 &  5.391e-04 &  8.277e-04 &  6.477e-04 &  1.867e+00 &  1.978e+00 &  1.842e+00 &  2.002e+00 &  1.801e+00 \\
	&  5.705e-03 &  1.507e-04 &  1.305e-04 &  1.462e-04 &  2.068e-04 &  1.772e-04 &  1.907e+00 &  1.961e+00 &  1.883e+00 &  2.001e+00 &  1.870e+00 \\
	&  2.853e-03 &  3.913e-05 &  3.326e-05 &  3.835e-05 &  5.166e-05 &  4.661e-05 &  1.946e+00 &  1.973e+00 &  1.931e+00 &  2.001e+00 &  1.927e+00 \\
	\hline
2	&  3.227e-02 &  2.533e-04 &  2.914e-04 &  2.886e-04 &  3.237e-04 &  2.072e-04 & - & - & - & - & - \\
	&  1.614e-02 &  2.999e-05 &  3.522e-05 &  3.911e-05 &  3.343e-05 &  2.972e-05 &  3.078e+00 &  3.048e+00 &  2.884e+00 &  3.275e+00 &  2.802e+00 \\
	&  8.069e-03 &  3.685e-06 &  4.334e-06 &  5.274e-06 &  3.494e-06 &  4.200e-06 &  3.025e+00 &  3.023e+00 &  2.891e+00 &  3.258e+00 &  2.823e+00 \\
	&  4.034e-03 &  4.636e-07 &  5.399e-07 &  6.953e-07 &  3.844e-07 &  5.662e-07 &  2.991e+00 &  3.005e+00 &  2.923e+00 &  3.184e+00 &  2.891e+00 \\
	&  2.017e-03 &  5.819e-08 &  6.739e-08 &  8.965e-08 &  4.365e-08 &  7.380e-08 &  \makegreen{2.994e+00} &  \makegreen{3.002e+00} &  \makegreen{2.955e+00} &  \makegreen{3.139e+00} &  \makegreen{2.940e+00} \\
	\hline
3	&  2.500e-02 &  3.565e-06 &  4.742e-06 &  4.378e-06 &  4.842e-06 &  2.831e-06 & - & - & - & - & - \\
	&  1.250e-02 &  2.183e-07 &  3.090e-07 &  2.854e-07 &  2.939e-07 &  1.867e-07 &  4.029e+00 &  3.940e+00 &  3.939e+00 &  4.042e+00 &  3.923e+00 \\
	&  6.250e-03 &  1.347e-08 &  1.948e-08 &  1.890e-08 &  1.731e-08 &  1.289e-08 &  4.018e+00 &  3.987e+00 &  3.917e+00 &  4.086e+00 &  3.856e+00 \\
	&  3.125e-03 &  7.961e-10 &  1.193e-09 &  1.232e-09 &  9.577e-10 &  8.673e-10 &  \makegreen{4.081e+00} &  \makegreen{4.030e+00} &  \makegreen{3.939e+00} &  \makegreen{4.176e+00} &  \makegreen{3.894e+00} \\
	\hline
\end{tabular}
}
\end{center}
\end{table}

From a previous study, it was noted that optimal convergence orders were not recovered for the $k1$ case with isoparametric geometry representation, even when the normals were defined according to the superparametric geometry representation~\cite[Section \makeblue{4.1.1}]{bassi1997}. We obtained similar results when using isoparametric geometry with exact normals, observing suboptimal convergence for the $k1$ case, but optimal convergence orders otherwise, even on the very high aspect ratio meshes. A thorough mathematical analysis is required to settle questions related to the need for superparametric geometry representation definitively.

\subsection{Results - Navier-Stokes}
\label{sec:verification_NavierStokes}

For the verification of the extension of the results obtained above to the Navier-Stokes equations, a 2D case with an analytical solution in an annular region employing a Neumann, no-slip adiabatic wall boundary condition was selected. It can be shown that an exact solution to the 2D Navier-Stokes equations is given by~\cite[Section \makeblue{7}]{illingworth1950}
\begin{align}
& T = T_i - \frac{C^2 \mu}{\kappa} \left( \frac{2}{ r_o^2} \log\left(\frac{r}{r_i}\right) + \left(\frac{1}{r^2}-\frac{1}{r_i^2}\right) \right),  \label{eq:TC_Analytical_T} \\
& \vect{v} = \begin{bmatrix} -\sin(\theta) & \cos(\theta) \end{bmatrix} v_{\theta}, \nonumber
\end{align}

where
\begin{align} \label{eq:TC_Analytical_v_theta}
C = \frac{\omega_i}{\frac{1}{r_i^2}-\frac{1}{r_o^2}},\ v_{\theta} = C r\left(\frac{1}{r^2}-\frac{1}{r_o^2} \right),\ \theta = \tan^{-1}\left(\frac{y}{x}\right),
\end{align}

and where $\mu$ is taken to be constant. In the interest of clarity, we provide a detailed derivation in~\autoref{sec:NS_2D_Analytical}. Above $i$ and $o$ subscripts are used to denote the inner and outer annulus surfaces. For the results computed below, we set $r_i = 0.5$, $r_o = 1.0$, $T_i = 1.0$, $\omega_i = 1.0$, $\rho_i = 1.0$, $Pr = 0.72$, $C_p = \frac{\gamma}{\gamma-1}$, $\gamma = 1.4$, and $\mu = 0.001$.

\begin{figure}
    \centering
    \begin{subfigure}[b]{0.485\textwidth}
        \includegraphics[width=\textwidth]{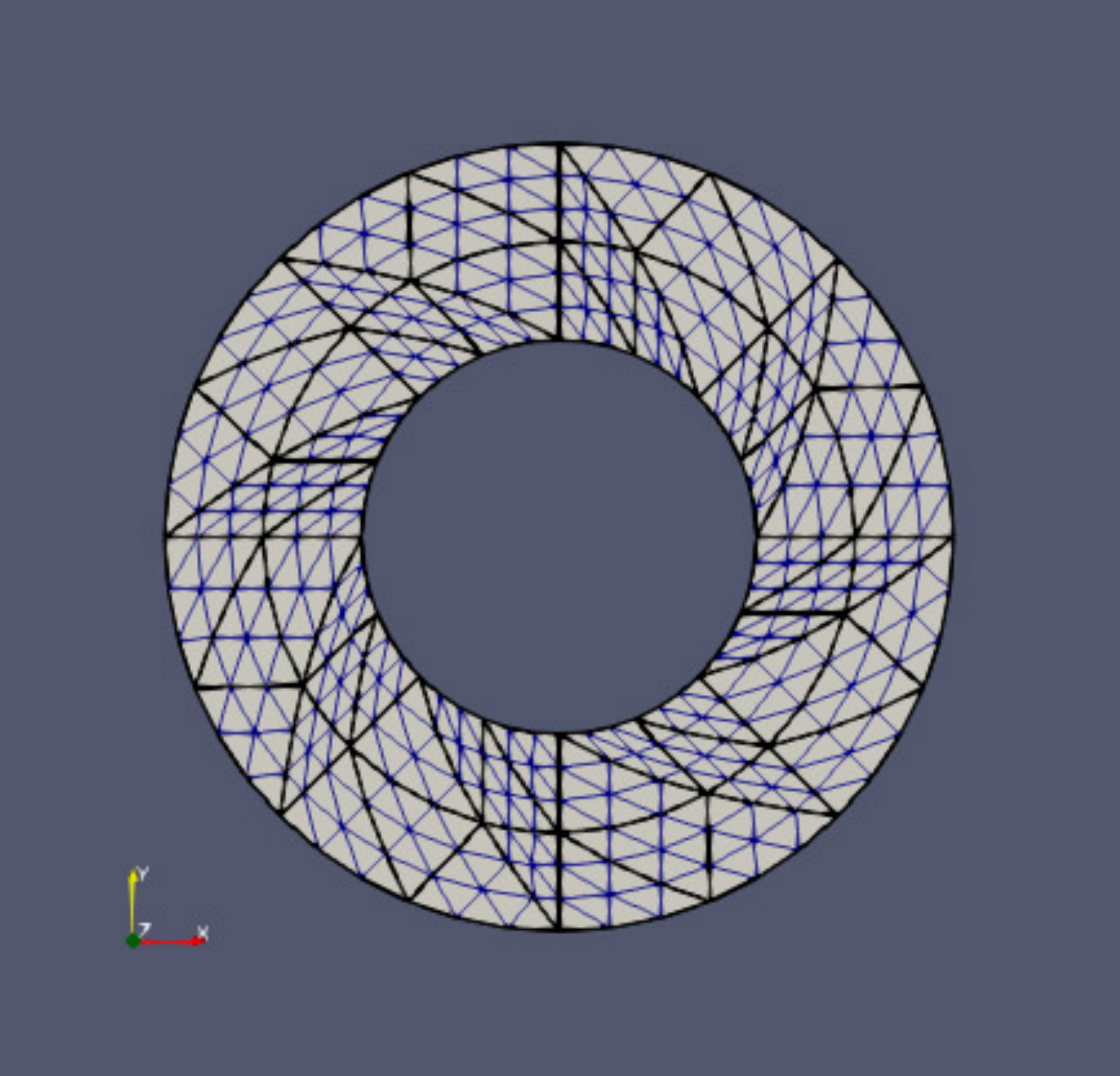}
        \caption{AR $ \approx 1.0$}
    \end{subfigure}
\hfill
    \centering
    \begin{subfigure}[b]{0.485\textwidth}
        \includegraphics[width=\textwidth]{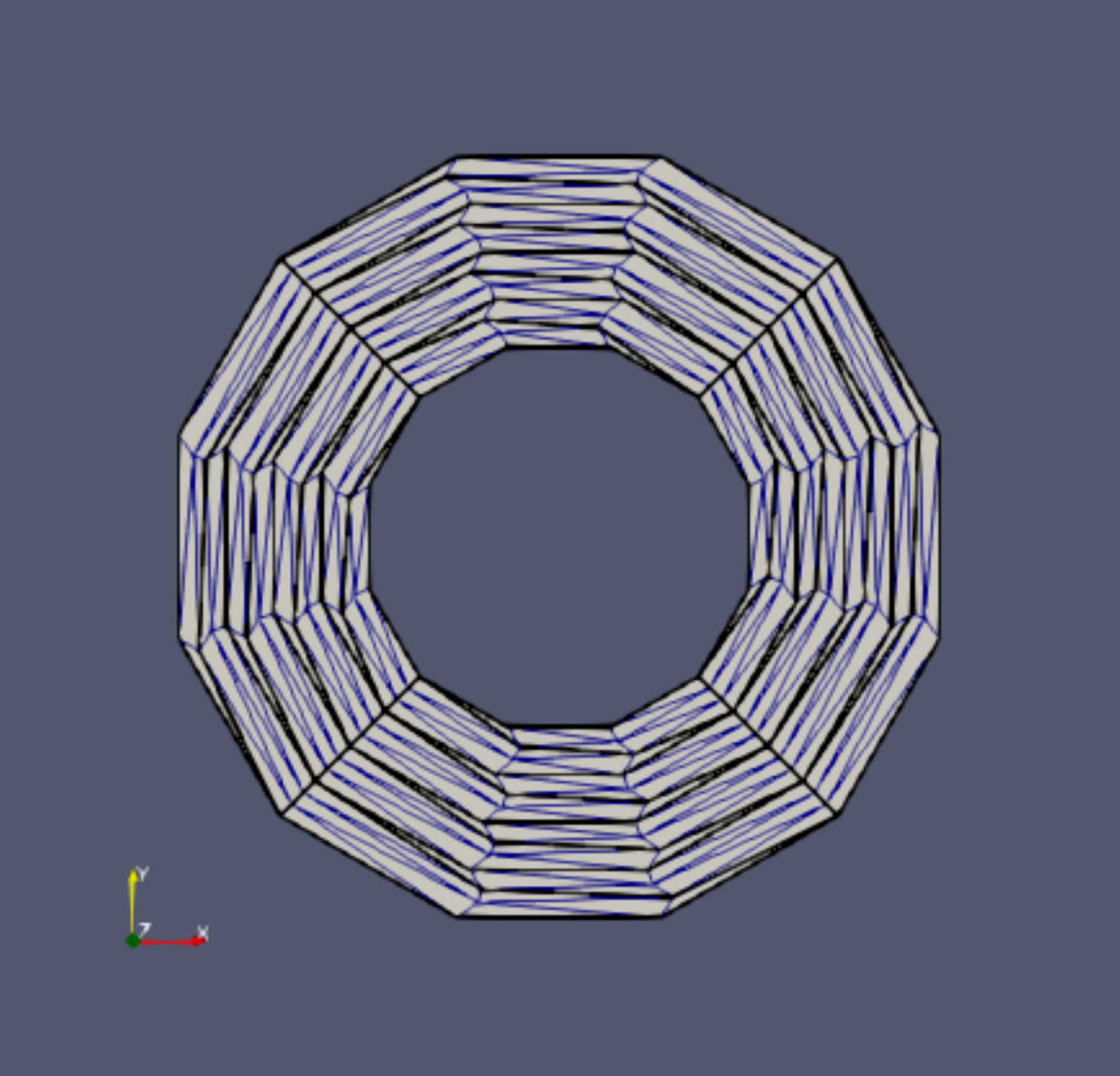}
        \caption{AR $ \approx 8.0$}
    \end{subfigure}
    \caption{Sample k3 Isoparametric Meshes used in the Navier-Stokes Convergence Studies}
    \label{fig:TaylorCouette_Meshes}
\end{figure}

No-slip isothermal and adiabatic boundary conditions were imposed on inner and outer walls, respectively. While it has been proven that the correct number of boundary conditions to impose for a no-slip wall to ensure well-posedness is equal to $d+1$ where $d$ is the dimension of the problem~\cite[Remark \makeblue{11} and Table \makeblue{2}]{nordstrom2005}, it can be noted in this case, that no boundary conditions are imposed for the remaining variable, leading to an ill-posed problem. Hence, a pressure Dirichlet boundary condition was additionally imposed on the inner wall. As the analytical solution is available, the integrated $L^2$ errors of the computed temperature and velocity components were used to assess optimal convergence. Similarly to the Euler case, we present the results on meshes comprised of elements having various aspect ratios with sample meshes shown in~\autoref{fig:TaylorCouette_Meshes}.

As expected from the theory established for isoparametric geometry representation for elliptic equations, we note that optimal solution convergence in $L^2$ was obtained for this case in the lower Reynolds number regime, when selecting $\mu = 1.0$, regardless of the mesh element aspect ratio. The results of~\autoref{tab:NS_TRI_mu1e-3_AR=1.0} and~\autoref{tab:NS_TRI_mu1e-3_AR=8.0} show that an isoparametric geometry representation remained sufficient to obtain optimal convergence orders when $\mu$ was decreased to $0.001$, however, in the case of the higher aspect ratio meshes, it can be noted that error magnitudes were significantly lower for several cases when using a superparametric geometry representation, most notably for the $k1$ and $k2$ solution orders (see~\autoref{tab:NS_TRI_mu1e-3_AR=8.0}).  It thus remains uncertain whether isoparametric geometry representation will continue to be sufficient for higher Reynolds number cases. We plan on performing this investigation in the future.

\begin{table}[!h]
    \centering
\begin{subtable}{.485\textwidth}
\begin{center}
\resizebox{\textwidth}{!}{
\begin{tabular}{| l | c | c c c | c c c | }
	\hline
	\multicolumn{2}{|c|}{} & \multicolumn{3}{c|}{$L^2$ Error} & \multicolumn{3}{c|}{Conv. Order} \\
	\hline
	Order ($k$) & Mesh Size ($h$) & $u$ & $v$ & $T$ & $u$ & $v$ & $T$ \\
	\hline
1	& 1.44e-01 & 2.78e-02 & 2.78e-02 & 1.47e-02 & -    & -    & -    \\
	& 7.22e-02 & 3.00e-02 & 3.00e-02 & 3.82e-03 & -0.11 & -0.11 & 1.94 \\
	& 3.61e-02 & 8.55e-03 & 8.55e-03 & 9.08e-04 & 1.81 & 1.81 & 2.07 \\
	& 1.80e-02 & 1.46e-03 & 1.46e-03 & 1.58e-04 & 2.55 & 2.55 & 2.52 \\
	& 9.02e-03 & 2.18e-04 & 2.18e-04 & 2.55e-05 & 2.75 & 2.75 & 2.63 \\
	\hline
2	& 1.02e-01 & 2.62e-02 & 2.62e-02 & 1.50e-03 & -    & -    & -    \\
	& 5.10e-02 & 3.73e-03 & 3.73e-03 & 3.34e-04 & 2.81 & 2.81 & 2.17 \\
	& 2.55e-02 & 2.53e-04 & 2.53e-04 & 2.73e-05 & 3.88 & 3.88 & 3.62 \\
	& 1.28e-02 & 1.76e-05 & 1.76e-05 & 1.80e-06 & 3.85 & 3.85 & 3.92 \\
	& 6.38e-03 & 1.79e-06 & 1.79e-06 & 1.76e-07 & 3.30 & 3.30 & 3.36 \\
	\hline
3	& 7.91e-02 & 3.90e-03 & 3.90e-03 & 3.30e-04 & -    & -    & -    \\
	& 3.95e-02 & 2.08e-04 & 2.08e-04 & 2.29e-05 & 4.23 & 4.23 & 3.85 \\
	& 1.98e-02 & 9.95e-06 & 9.95e-06 & 1.17e-06 & 4.38 & 4.38 & 4.29 \\
	& 9.88e-03 & 5.84e-07 & 5.84e-07 & 7.76e-08 & 4.09 & 4.09 & 3.91 \\
	& 4.94e-03 & 3.53e-08 & 3.53e-08 & 5.11e-09 & 4.05 & 4.05 & 3.93 \\
	\hline
\end{tabular}
}
\end{center}
\caption{$k_G = k$}
\end{subtable}%
\hfill
\begin{subtable}{.485\textwidth}
\begin{center}
\resizebox{\textwidth}{!}{
\begin{tabular}{| l | c | c c c | c c c | }
	\hline
	\multicolumn{2}{|c|}{} & \multicolumn{3}{c|}{$L^2$ Error} & \multicolumn{3}{c|}{Conv. Order} \\
	\hline
	Order ($k$) & Mesh Size ($h$) & $u$ & $v$ & $T$ & $u$ & $v$ & $T$ \\
	\hline
1	& 1.44e-01 & 4.34e-02 & 4.34e-02 & 4.06e-03 & -    & -    & -    \\
	& 7.22e-02 & 2.37e-02 & 2.37e-02 & 1.61e-03 & 0.87 & 0.87 & 1.34 \\
	& 3.61e-02 & 8.13e-03 & 8.13e-03 & 6.14e-04 & 1.54 & 1.54 & 1.39 \\
	& 1.80e-02 & 1.56e-03 & 1.56e-03 & 1.35e-04 & 2.38 & 2.38 & 2.19 \\
	& 9.02e-03 & 2.70e-04 & 2.70e-04 & 2.65e-05 & 2.52 & 2.52 & 2.35 \\
	\hline
2	& 1.02e-01 & 2.85e-02 & 2.85e-02 & 1.80e-03 & -    & -    & -    \\
	& 5.10e-02 & 3.86e-03 & 3.86e-03 & 3.56e-04 & 2.88 & 2.88 & 2.34 \\
	& 2.55e-02 & 2.48e-04 & 2.48e-04 & 2.75e-05 & 3.96 & 3.96 & 3.69 \\
	& 1.28e-02 & 1.58e-05 & 1.58e-05 & 1.61e-06 & 3.97 & 3.97 & 4.09 \\
	& 6.38e-03 & 1.55e-06 & 1.55e-06 & 1.43e-07 & 3.36 & 3.36 & 3.49 \\
	\hline
3	& 7.91e-02 & 4.04e-03 & 4.04e-03 & 3.58e-04 & -    & -    & -    \\
	& 3.95e-02 & 2.05e-04 & 2.05e-04 & 2.30e-05 & 4.30 & 4.30 & 3.96 \\
	& 1.98e-02 & 9.09e-06 & 9.09e-06 & 1.07e-06 & 4.49 & 4.49 & 4.43 \\
	& 9.88e-03 & 5.20e-07 & 5.20e-07 & 6.90e-08 & 4.13 & 4.13 & 3.95 \\
	& 4.94e-03 & 3.13e-08 & 3.13e-08 & 4.52e-09 & 4.06 & 4.06 & 3.93 \\
	\hline
\end{tabular}
}
\end{center}
\caption{$k_G = k+1$}
\end{subtable}
\caption{Errors and Convergence Orders - ToBeCurvedTRI Meshes ($\mu = 0.001$, AR $\approx$ 1.0)}
\label{tab:NS_TRI_mu1e-3_AR=1.0}
\end{table}

\begin{table}[!h]
    \centering
\begin{subtable}{.485\textwidth}
\begin{center}
\resizebox{\textwidth}{!}{
\begin{tabular}{| l | c | c c c | c c c | }
	\hline
	\multicolumn{2}{|c|}{} & \multicolumn{3}{c|}{$L^2$ Error} & \multicolumn{3}{c|}{Conv. Order} \\
	\hline
	Order ($k$) & Mesh Size ($h$) & $u$ & $v$ & $T$ & $u$ & $v$ & $T$ \\
	\hline
1	& 1.02e-01 & 1.38e-01 & 1.38e-01 & 9.01e-02 & -    & -    & -    \\
	& 5.10e-02 & 4.29e-02 & 4.29e-02 & 2.34e-02 & 1.69 & 1.69 & 1.94 \\
	& 2.55e-02 & 1.46e-02 & 1.46e-02 & 4.94e-03 & 1.56 & 1.56 & 2.25 \\
	& 1.28e-02 & \makered{3.29e-03} & \makered{3.29e-03} & \makered{7.07e-04} & 2.15 & 2.15 & 2.80 \\
	\hline
2	& 7.22e-02 & 5.55e-03 & 5.55e-03 & 1.00e-03 & -    & -    & -    \\
	& 3.61e-02 & 4.65e-04 & 4.65e-04 & 1.07e-04 & 3.58 & 3.58 & 3.23 \\
	& 1.80e-02 & 4.43e-05 & 4.43e-05 & 1.16e-05 & 3.39 & 3.39 & 3.20 \\
	& 9.02e-03 & \makered{5.37e-06} & \makered{5.37e-06} & \makered{1.19e-06} & 3.04 & 3.04 & 3.29 \\
	\hline
3	& 5.59e-02 & 8.80e-04 & 8.80e-04 & 1.10e-04 & -    & -    & -    \\
	& 2.79e-02 & 5.35e-05 & 5.35e-05 & 6.48e-06 & 4.04 & 4.04 & 4.09 \\
	& 1.40e-02 & 3.40e-06 & 3.40e-06 & 4.11e-07 & 3.98 & 3.98 & 3.98 \\
	& 6.99e-03 & 2.19e-07 & 2.19e-07 & 2.59e-08 & 3.96 & 3.96 & 3.99 \\
	\hline
\end{tabular}
}
\end{center}
\caption{$k_G = k$}
\end{subtable}%
\hfill
\begin{subtable}{.485\textwidth}
\begin{center}
\resizebox{\textwidth}{!}{
\begin{tabular}{| l | c | c c c | c c c | }
	\hline
	\multicolumn{2}{|c|}{} & \multicolumn{3}{c|}{$L^2$ Error} & \multicolumn{3}{c|}{Conv. Order} \\
	\hline
	Order ($k$) & Mesh Size ($h$) & $u$ & $v$ & $T$ & $u$ & $v$ & $T$ \\
	\hline
1	& 1.02e-01 & 4.72e-02 & 4.72e-02 & 3.39e-03 & -    & -    & -    \\
	& 5.10e-02 & 1.34e-02 & 1.34e-02 & 9.63e-04 & 1.81 & 1.81 & 1.82 \\
	& 2.55e-02 & 2.57e-03 & 2.57e-03 & 2.31e-04 & 2.39 & 2.39 & 2.06 \\
	& 1.28e-02 & \makegreen{4.06e-04} & \makegreen{4.06e-04} & \makegreen{4.97e-05} & 2.66 & 2.66 & 2.22 \\
	\hline
2	& 7.22e-02 & 5.38e-03 & 5.38e-03 & 6.05e-04 & -    & -    & -    \\
	& 3.61e-02 & 3.65e-04 & 3.65e-04 & 6.15e-05 & 3.88 & 3.88 & 3.30 \\
	& 1.80e-02 & 3.75e-05 & 3.75e-05 & 7.21e-06 & 3.28 & 3.28 & 3.09 \\
	& 9.02e-03 & \makegreen{4.70e-06} & \makegreen{4.70e-06} & \makegreen{8.12e-07} & 3.00 & 3.00 & 3.15 \\
	\hline
3	& 5.59e-02 & 7.87e-04 & 7.87e-04 & 1.34e-04 & -    & -    & -    \\
	& 2.79e-02 & 4.61e-05 & 4.61e-05 & 7.15e-06 & 4.09 & 4.09 & 4.22 \\
	& 1.40e-02 & 2.98e-06 & 2.98e-06 & 4.25e-07 & 3.95 & 3.95 & 4.07 \\
	& 6.99e-03 & 1.99e-07 & 1.99e-07 & 2.52e-08 & 3.91 & 3.91 & 4.08 \\
	\hline
\end{tabular}
}
\end{center}
\caption{$k_G = k+1$}
\end{subtable}
\caption{Errors and Convergence Orders - ToBeCurvedTRI Meshes ($\mu = 0.001$, AR $\approx$ 8.0)}
\label{tab:NS_TRI_mu1e-3_AR=8.0}
\end{table}

The analogous results obtained when using quadrilateral elements are presented in~\autoref{sec:NS_QUAD_results}.

\section{Conclusion}
\label{sec:conclusion}

The need for superparametric geometry representation when solving the Euler equations with slip wall boundary conditions has been demonstrated numerically. Despite not requiring a superparametric geometry representation when solving the Navier-Stokes equations, the results obtained for the case presented here demonstrated that significantly lower errors may be obtained when employing the higher-order geometry representation, motivating a follow-up investigation in the higher Reynolds number regime. Finally, we note that the additional cost associated with the increased geometry representation order was negligible in all cases. Extension to the 3D case will be pursued in future work.

\section*{Acknowledgements}

The authors gratefully acknowledge the support from the Natural Sciences and Engineering Research Council (NSERC) and the Department of Mechanical Engineering of McGill University.

\begin{appendices}

\section{Navier-Stokes 2D Analytical Solution - Annular Domain}
\label{sec:NS_2D_Analytical}

First, we express the 2D steady continuity, Navier-Stokes and Energy equations in polar coordinates. The continuity equation is given by~\cite[eq. (\makeblue{3.70})]{langlois2014}
\begin{align} \label{eq:2D_NS_cyl_continuity}
& v_r \frac{\partial \rho}{\partial r} + \frac{v_{\theta}}{r} \frac{\partial \rho}{\partial \theta} + \rho \Delta = 0,
\end{align}

where the dilation, $\Delta$, is given by
\begin{align*}
& \Delta = \frac{\partial v_r}{\partial r} + \frac{v_r}{r} + \frac{1}{r} \frac{\partial v_{\theta}}{\partial \theta}.
\end{align*}

For the Navier-Stokes equations, we have~\cite[eq. (\makeblue{3.72})]{langlois2014}
\begin{alignat}{3}
& \rho \left( v_r \frac{\partial v_r}{\partial r} + \frac{v_{\theta}}{r} \frac{\partial v_r}{\partial \theta} -\frac{v_{\theta}^2}{r} \right) 
&& = \rho f_r - \frac{\partial p}{\partial r} 
+ (\lambda + \mu)  \frac{\partial \Delta}{\partial r} + \mu \left( \oldnabla^2 v_r - \frac{v_r}{r^2} - \frac{2}{r^2}  \frac{\partial v_{\theta}}{\partial \theta}  \right) \nonumber \\
&&& + \Delta \frac{\partial \lambda}{\partial r} + 2 \left( e_{rr} \frac{\partial \mu}{\partial r} + \frac{1}{r} e_{r\theta}  \frac{\partial \mu}{\partial \theta} \right)
\label{eq:2D_NS_cyl_mom1} \\ 
& \rho \left( v_r \frac{\partial v_{\theta}}{\partial r} + \frac{v_{\theta}}{r} \frac{\partial v_{\theta}}{\partial \theta} +\frac{v_r v_{\theta}}{r} \right) 
&& = \rho f_{\theta} - \frac{1}{r} \frac{\partial p}{\partial {\theta}} 
+ \frac{\lambda + \mu}{r}  \frac{\partial \Delta}{\partial {\theta}} + \mu \left( \oldnabla^2 v_{\theta} + \frac{2}{r^2}  \frac{\partial v_r}{\partial \theta}  - \frac{v_{\theta}}{r^2} \right) \nonumber \\
&&& + \frac{\Delta}{r} \frac{\partial \lambda}{\partial {\theta}} + 2 \left( e_{{\theta}r} \frac{\partial \mu}{\partial r} + \frac{1}{r} e_{{\theta}\theta}  \frac{\partial \mu}{\partial \theta} \right) \label{eq:2D_NS_cyl_mom2}
\end{alignat}

where
\begin{align*}
& \oldnabla^2 = \frac{\partial^2}{\partial r^2} + \frac{1}{r} \frac{\partial}{\partial r} + \frac{1}{r^2} \frac{\partial^2}{\partial \theta^2},   
\end{align*}

and the rate of deformation tensor components are defined according to
\begin{align*}
& e_{rr} = \frac{\partial v_r}{\partial r},  \\
& e_{\theta\theta} = \frac{1}{r} \left(\frac{\partial v_{\theta}}{\partial \theta} + v_r \right),  \\
& e_{r\theta} = e_{\theta r}  = \frac{1}{2} \left( \frac{\partial v_{\theta}}{\partial r} - \frac{v_{\theta}}{r} + \frac{1}{r} \frac{\partial v_r}{\partial \theta} \right).
\end{align*}

Finally, the energy equation is given by
\begin{align} \label{eq:2D_NS_cyl_ener}
& \rho \left( v_r \frac{\partial e}{\partial r} + \frac{v_{\theta}}{r} \frac{\partial e}{\partial \theta} \right) 
= - p \Delta
+ \frac{\partial}{\partial r} \left( \kappa \frac{\partial T}{\partial r} \right) + \frac{\kappa}{r} \frac{\partial T}{\partial r}
+ \frac{1}{r^2} \frac{\partial}{\partial \theta} \left( \kappa \frac{\partial T}{\partial \theta} \right) + \lambda \Delta^2 + 2\mu \sum_{\alpha,\beta = r,\theta} e_{\alpha\beta}^2,
\end{align}


where $e$ is the internal energy. Making the following assumptions:

\begin{itemize}
\item the radial component of the velocity is zero, $v_r = 0$;
\item the solution is radially symmetric, $\frac{\partial}{\partial \theta} ( \cdot ) = 0$;
\item the body forces are negligible;
\item the viscosity is constant, $\mu = \text{const.}$;
\item the coefficient of thermal conductivity is constant, $\kappa = \text{const.}$;
\item the gas is calorically ideal,
\end{itemize}

it can first be noted that the dilation is zero, reducing~\eqref{eq:2D_NS_cyl_continuity} to
\begin{align*}
0 = 0,
\end{align*}

implying that the continuity equation is satisfied for all choices of $\rho$. Similarly,~\eqref{eq:2D_NS_cyl_mom1} is simplified to
\begin{align} \label{eq:2D_NS_cyl_mom1_simp}
\frac{\rho v_{\theta}^2}{r}
= \frac{d p}{d r} ,
\end{align}

and~\eqref{eq:2D_NS_cyl_mom2} is simplified to
\begin{align} \label{eq:2D_NS_cyl_mom2_simp}
& \frac{d^2 v_{\theta}}{d r^2} + \frac{1}{r} \frac{d v_{\theta}}{d r} - \frac{v_{\theta}}{r^2} = 0 \nonumber \\ 
\rightarrow\ & r^2 \frac{d^2 v_{\theta}}{d r^2} + r \frac{d v_{\theta}}{d r} - v_{\theta} = \frac{d}{d r} \left[r^2 \left(\frac{d v_{\theta}}{d r} - \frac{v_{\theta}}{r} \right) \right] = 0.
\end{align}

Finally, the energy equation, \eqref{eq:2D_NS_cyl_ener}, reduces to
\begin{align} \label{eq:2D_NS_cyl_ener_simp}
& \mu \left( \frac{1}{Pr} \frac{d^2 h}{d r^2} + \frac{1}{r Pr} \frac{d h}{d r} + \left( \frac{d v_{\theta}}{d r} - \frac{v_{\theta}}{r} \right)^2 \right) = 0 \nonumber \\
\rightarrow\ &
\frac{1}{r Pr} \frac{d}{d r} \left( r \frac{d h}{d r} \right) + \left( \frac{d v_{\theta}}{d r} - \frac{v_{\theta}}{r} \right)^2 
= 0,
\end{align}

where
\begin{align*}
h = C_p T.
\end{align*}

Defining $\omega \coloneqq \frac{v_{\theta}}{r}$, it can be noted that
\begin{align} \label{eq:identity_omega}
\frac{d v_{\theta}}{d r} - \frac{v_{\theta}}{r} = \omega + r \frac{d \omega}{d r} - \omega = r \frac{d \omega}{d r},
\end{align}

such that~\eqref{eq:2D_NS_cyl_mom2_simp} can then be expressed as
\begin{align} \label{eq:TC_const1}
r^3 \frac{d \omega}{d r} = \text{const.} \coloneqq C_1.
\end{align}

Using~\eqref{eq:identity_omega} and~\eqref{eq:TC_const1},~\eqref{eq:2D_NS_cyl_ener_simp} can be rewritten as
\begin{align*}
\frac{d}{d r} \left( r \frac{d h}{d r} \right) 
= -\frac{Pr}{r^3}  \left( r^3 \frac{d \omega}{d r} \right)^2 
= -\frac{Pr}{r^3}  C_1^2.
\end{align*}

Integrating with respect to $r$, we obtain
\begin{align*}
r \frac{d h}{d r} -\frac{Pr C_1^2}{2 r^2} = r \left( \frac{d h}{d r} + r^2 Pr  \frac{d \omega}{d r} \left(-\frac{r}{2} \frac{d \omega}{d r} \right) \right) = \text{const.}
\end{align*}

Assuming that 
\begin{align} \label{eq:TC_Assumption_omega}
-\frac{r}{2} \frac{d \omega}{d r} = \omega + \text{const.}, 
\end{align}

and noting~\eqref{eq:TC_const1},the last condition from Illingworth is recovered~\cite[eq. (\makeblue{47})]{illingworth1950},
\begin{align} \label{eq:TC_const2}
r \left( \frac{d h}{d r} + r^2 Pr  \omega \frac{d \omega}{d r} \right) = \text{const.}
\end{align}

With the assumptions made above, it is then straightforward to obtain the solution of~\autoref{sec:verification_NavierStokes}. Noting the constraint of~\eqref{eq:TC_const1}, the angular velocity may take the form,
\begin{align*}
\omega = \omega_i \frac{\frac{1}{r^2} + C_{\omega}}{\frac{1}{r_i^2} + C_{\omega}},
\end{align*}

which also satisfies~\eqref{eq:TC_Assumption_omega}. Setting the two no-slip boundary conditions with prescribed non-zero angular velocity, $\omega_i$, at $r = r_i$ and zero angular velocity at $r = r_o$ gives the expression for $v_{\theta}$ in~\eqref{eq:TC_Analytical_v_theta}. Substituting into~\eqref{eq:TC_const2} and imposing the Dirichlet temperature boundary condition on the inner wall and the Neumann adiabatic boundary condition on the outer wall then leads to the solution for the temperature distribution given in~\eqref{eq:TC_Analytical_T}.

\section{Navier-Stokes Results on Quadrilateral Meshes}
\label{sec:NS_QUAD_results}

In this appendix, we present results of the convergence order study for the Navier-Stokes test case of~\autoref{sec:verification_NavierStokes} computed using meshes consisting of quadrilateral elements. As for the results on the triangular meshes, similar trends were observed when selecting $\mu = 1.0$.

\begin{table}[!ht]
    \centering
\begin{subtable}{.485\textwidth}
\begin{center}
\resizebox{\textwidth}{!}{
\begin{tabular}{| l | c | c c c | c c c | }
	\hline
	\multicolumn{2}{|c|}{} & \multicolumn{3}{c|}{$L^2$ Error} & \multicolumn{3}{c|}{Conv. Order} \\
	\hline
	Order ($k$) & Mesh Size ($h$) & $u$ & $v$ & $T$ & $u$ & $v$ & $T$ \\
	\hline
1	& 1.77e-01 & 2.25e-02 & 2.25e-02 & 3.00e-02 & -    & -    & -    \\
	& 8.84e-02 & 7.09e-03 & 7.09e-03 & 6.97e-03 & 1.67 & 1.67 & 2.11 \\
	& 4.42e-02 & 1.13e-03 & 1.13e-03 & 1.06e-03 & 2.64 & 2.64 & 2.72 \\
	& 2.21e-02 & 2.38e-04 & 2.38e-04 & 1.46e-04 & 2.26 & 2.26 & 2.87 \\
	& 1.11e-02 & 5.50e-05 & 5.50e-05 & 2.30e-05 & 2.11 & 2.11 & 2.66 \\
	\hline
2	& 1.18e-01 & 2.83e-02 & 2.83e-02 & 2.04e-03 & -    & -    & -    \\
	& 5.89e-02 & 1.66e-03 & 1.66e-03 & 1.76e-04 & 4.09 & 4.09 & 3.54 \\
	& 2.95e-02 & 8.12e-05 & 8.12e-05 & 1.15e-05 & 4.35 & 4.35 & 3.93 \\
	& 1.47e-02 & 5.82e-06 & 5.82e-06 & 9.42e-07 & 3.80 & 3.80 & 3.62 \\
	& 7.37e-03 & 6.58e-07 & 6.58e-07 & 1.10e-07 & 3.15 & 3.15 & 3.10 \\
	\hline
3	& 8.84e-02 & 3.14e-03 & 3.14e-03 & 2.99e-04 & -    & -    & -    \\
	& 4.42e-02 & 8.33e-05 & 8.33e-05 & 1.07e-05 & 5.24 & 5.24 & 4.80 \\
	& 2.21e-02 & 2.19e-06 & 2.19e-06 & 3.20e-07 & 5.25 & 5.25 & 5.07 \\
	& 1.11e-02 & 9.61e-08 & 9.61e-08 & 1.65e-08 & 4.51 & 4.51 & 4.27 \\
	& 5.52e-03 & 5.06e-09 & 5.06e-09 & 1.04e-09 & 4.25 & 4.25 & 3.99 \\
	\hline
\end{tabular}
}
\end{center}
\caption{$k_G = k$}
\end{subtable}%
\hfill
\begin{subtable}{.485\textwidth}
\begin{center}
\resizebox{\textwidth}{!}{
\begin{tabular}{| l | c | c c c | c c c | }
	\hline
	\multicolumn{2}{|c|}{} & \multicolumn{3}{c|}{$L^2$ Error} & \multicolumn{3}{c|}{Conv. Order} \\
	\hline
	Order ($k$) & Mesh Size ($h$) & $u$ & $v$ & $T$ & $u$ & $v$ & $T$ \\
	\hline
1	& 1.77e-01 & 2.23e-02 & 2.23e-02 & 3.12e-03 & -    & -    & -    \\
	& 8.84e-02 & 1.04e-02 & 1.04e-02 & 9.15e-04 & 1.11 & 1.11 & 1.77 \\
	& 4.42e-02 & 2.98e-03 & 2.98e-03 & 3.26e-04 & 1.80 & 1.80 & 1.49 \\
	& 2.21e-02 & 6.21e-04 & 6.21e-04 & 8.28e-05 & 2.26 & 2.26 & 1.98 \\
	& 1.11e-02 & 1.33e-04 & 1.33e-04 & 2.02e-05 & 2.23 & 2.23 & 2.03 \\
	\hline
2	& 1.18e-01 & 2.84e-02 & 2.84e-02 & 2.05e-03 & -    & -    & -    \\
	& 5.89e-02 & 1.67e-03 & 1.67e-03 & 1.76e-04 & 4.09 & 4.09 & 3.54 \\
	& 2.95e-02 & 8.20e-05 & 8.20e-05 & 1.16e-05 & 4.35 & 4.35 & 3.93 \\
	& 1.47e-02 & 5.90e-06 & 5.90e-06 & 9.43e-07 & 3.80 & 3.80 & 3.62 \\
	& 7.37e-03 & 6.66e-07 & 6.66e-07 & 1.10e-07 & 3.15 & 3.15 & 3.10 \\
	\hline
3	& 8.84e-02 & 3.14e-03 & 3.14e-03 & 2.98e-04 & -    & -    & -    \\
	& 4.42e-02 & 8.32e-05 & 8.32e-05 & 1.07e-05 & 5.24 & 5.24 & 4.80 \\
	& 2.21e-02 & 2.19e-06 & 2.19e-06 & 3.19e-07 & 5.25 & 5.25 & 5.07 \\
	& 1.11e-02 & 9.60e-08 & 9.60e-08 & 1.65e-08 & 4.51 & 4.51 & 4.27 \\
	& 5.52e-03 & 5.06e-09 & 5.06e-09 & 1.04e-09 & 4.25 & 4.25 & 3.98 \\
	\hline
\end{tabular}
}
\end{center}
\caption{$k_G = k+1$}
\end{subtable}
\caption{Errors and Convergence Orders - ToBeCurvedQUAD Meshes ($\mu = 0.001$, AR $\approx$ 1.0)}
\label{tab:NS_QUAD_mu1e-3_AR=1.0}
\end{table}

\begin{table}[!ht]
    \centering
\begin{subtable}{.485\textwidth}
\begin{center}
\resizebox{\textwidth}{!}{
\begin{tabular}{| l | c | c c c | c c c | }
	\hline
	\multicolumn{2}{|c|}{} & \multicolumn{3}{c|}{$L^2$ Error} & \multicolumn{3}{c|}{Conv. Order} \\
	\hline
	Order ($k$) & Mesh Size ($h$) & $u$ & $v$ & $T$ & $u$ & $v$ & $T$ \\
	\hline
1	& 1.25e-01 & 1.50e-01 & 1.50e-01 & 8.89e-02 & -    & -    & -    \\
	& 6.25e-02 & 5.89e-02 & 5.89e-02 & 3.33e-02 & 1.35 & 1.35 & 1.42 \\
	& 3.12e-02 & 1.95e-02 & 1.95e-02 & 7.27e-03 & 1.59 & 1.59 & 2.19 \\
	& 1.56e-02 & 4.14e-03 & 4.14e-03 & 9.90e-04 & 2.24 & 2.24 & 2.88 \\
	\hline
2	& 8.33e-02 & 5.80e-03 & 5.80e-03 & 6.37e-04 & -    & -    & -    \\
	& 4.17e-02 & 3.84e-04 & 3.84e-04 & 4.33e-05 & 3.92 & 3.92 & 3.88 \\
	& 2.08e-02 & 2.44e-05 & 2.44e-05 & 2.85e-06 & 3.98 & 3.98 & 3.92 \\
	& 1.04e-02 & 1.53e-06 & 1.53e-06 & 1.84e-07 & 3.99 & 3.99 & 3.95 \\
	\hline
3	& 6.25e-02 & 7.87e-04 & 7.87e-04 & 2.17e-04 & -    & -    & -    \\
	& 3.12e-02 & 5.04e-05 & 5.04e-05 & 1.22e-05 & 3.96 & 3.96 & 4.15 \\
	& 1.56e-02 & 2.39e-06 & 2.39e-06 & 5.73e-07 & 4.40 & 4.40 & 4.42 \\
	& 7.81e-03 & 9.87e-08 & 9.87e-08 & 2.32e-08 & 4.59 & 4.59 & 4.63 \\
	\hline
\end{tabular}
}
\end{center}
\caption{$k_G = k$}
\end{subtable}%
\hfill
\begin{subtable}{.485\textwidth}
\begin{center}
\resizebox{\textwidth}{!}{
\begin{tabular}{| l | c | c c c | c c c | }
	\hline
	\multicolumn{2}{|c|}{} & \multicolumn{3}{c|}{$L^2$ Error} & \multicolumn{3}{c|}{Conv. Order} \\
	\hline
	Order ($k$) & Mesh Size ($h$) & $u$ & $v$ & $T$ & $u$ & $v$ & $T$ \\
	\hline
1	& 1.25e-01 & 7.51e-02 & 7.51e-02 & 3.79e-03 & -    & -    & -    \\
	& 6.25e-02 & 2.49e-02 & 2.49e-02 & 1.57e-03 & 1.59 & 1.59 & 1.28 \\
	& 3.12e-02 & 4.21e-03 & 4.21e-03 & 3.54e-04 & 2.56 & 2.56 & 2.15 \\
	& 1.56e-02 & 5.69e-04 & 5.69e-04 & 6.81e-05 & 2.89 & 2.89 & 2.38 \\
	\hline
2	& 8.33e-02 & 1.73e-03 & 1.73e-03 & 1.41e-04 & -    & -    & -    \\
	& 4.17e-02 & 1.92e-04 & 1.92e-04 & 1.97e-05 & 3.17 & 3.17 & 2.84 \\
	& 2.08e-02 & 3.25e-05 & 3.25e-05 & 2.44e-06 & 2.56 & 2.56 & 3.01 \\
	& 1.04e-02 & 5.09e-06 & 5.09e-06 & 2.31e-07 & 2.68 & 2.68 & 3.40 \\
	\hline
3	& 6.25e-02 & 1.58e-04 & 1.58e-04 & 3.75e-05 & -    & -    & -    \\
	& 3.12e-02 & 7.60e-06 & 7.60e-06 & 2.51e-06 & 4.37 & 4.37 & 3.90 \\
	& 1.56e-02 & 3.63e-07 & 3.63e-07 & 1.27e-07 & 4.39 & 4.39 & 4.30 \\
	& 7.81e-03 & 1.75e-08 & 1.75e-08 & 5.31e-09 & 4.38 & 4.38 & 4.58 \\
	\hline
\end{tabular}
}
\end{center}
\caption{$k_G = k+1$}
\end{subtable}
\caption{Errors and Convergence Orders - ToBeCurvedQUAD Meshes ($\mu = 0.001$, AR $\approx$ 8.0)}
\label{tab:NS_QUAD_mu1e-3_AR=8.0}
\end{table}

\end{appendices}

\bibliography{Global}
\bibliographystyle{aiaa}

\end{document}